\documentclass[12pt]{article}
\usepackage[russian]{babel}
\usepackage{amsthm,amsfonts,amssymb,amscd}

\begin{document}
\begin{center}
\Large \bf Birationally rigid varieties \\
with a pencil of Fano double covers. I
\end{center}
\vspace{0.7cm}

\centerline{\large \bf Aleksandr V. Pukhlikov}%\footnote{ (
% 02-15-99258)}
\vspace{0.7cm}

\begin{center}
Max-Planck-Institut f\" ur Mathematik \\
Vivatsgasse 7 \\
53111 Bonn \\
GERMANY \\
e-mail: {\it pukh@mpim-bonn.mpg.de}
\end{center}
\vspace{0.3cm}

\begin{center}
Steklov Institute of Mathematics \\
Gubkina 8 \\
117966 Moscow \\
RUSSIA \\
e-mail: {\it pukh@mi.ras.ru}
\end{center}
\vspace{0.3cm}

\begin{center}
Division of Pure Mathematics \\
Department of Mathematical Sciences \\
M$\&$O Building, Peach Street \\
The University of Liverpool \\
Liverpool L69 7ZL \\
ENGLAND\\
e-mail: {\it pukh@liv.ac.uk}
\end{center}
\vspace{1cm}

\centerline{October 3, 2003}\vspace{1cm}

\parshape=1
3cm 10cm \noindent {\small \quad \quad \quad
\quad\quad\quad\quad\quad\quad\quad {\bf Abstract}\newline\newline
We prove that a general Fano fibration $\pi\colon V\to {\mathbb
P}^1$, the fiber of which is a double Fano hypersurface of index
1, is birationally superrigid provided it is sufficiently twisted
over the base. In particular, on $V$ there are no other structures
of a rationally connected fibration. The proof is based on the
method of maximal singularities.} \vspace{1cm}

\newpage

CONTENTS \vspace{0.5cm}

\noindent Introduction

0.1. Birationally rigid varieties

0.2. Varieties with a pencil of double covers

0.3. The main result

0.4. Historical remarks

0.5. Acknowledgements \vspace{0.3cm}

\noindent 1. The method of maximal singularities and the
regularity conditions

1.1. A criterion of birational rigidity

1.2. An explicit construction of the fibration $V/{\mathbb P}^1$

1.3. The regularity conditions outside the branch divisor

1.4. The regularity conditions on the branch divisor

1.5. Start of the proof of Theorem 1 \vspace{0.3cm}

\noindent 2. Singularity of a fiber outside the branch divisor

2.1. Hypertangent divisors and linear systems

2.2. Scheme of the proof of the condition (vs)

2.3. Movable families of curves \vspace{0.3cm}

\noindent 3. Singularity of a fiber on the branch divisor

3.1. Notations and discussion of the regularity condition

3.2. Start of the proof of the condition (vs)

3.3. Hypertangent divisors and tangent cones

3.4. Constructing new cycles

3.5. Degrees and multiplicities \vspace{0.3cm}

\noindent References

\newpage

\section*{Introduction}

In this paper we study birational geometry of higher-dimensional
algebraic varieties with a pencil of Fano double covers. The main
result of the paper, that is, the theorem on birational
superrigidity of these varieties provided they are sufficiently
twisted over the base, is formulated below in Sec. 0.3. The paper
presents the outcome of the first stage of the research; the
second part of this work, dealing with a relaxation of the
twistedness condition, will be published in the subsequent paper.

\subsection{Birationally rigid varieties}

A rationally connected projective variety $V$ with at most
${\mathbb Q}$-factorial terminal singularities is said to be {\it
birationally rigid}, if for any birational map
$$
\chi\colon V - \,-\,\to V',
$$
where $V'$ belongs to the same class of varieties, and any moving
linear system $\Sigma'$ on $V'$ there exists a birational self-map
$\chi^*\in \mathop{\rm Bir}V$, providing the following inequality
\begin{equation}
\label{0.1} c(\Sigma)\leq c(\Sigma'),
\end{equation}
where $\Sigma=(\chi\circ\chi^*)_*\Sigma'$ is the strict transform
of the linear system $\Sigma'$ on $V$ with respect to the
birational map
$$
\chi\circ\chi^*\colon V\stackrel{\chi^*}{- \,-\,\to}
V\stackrel{\chi}{- \,-\,\to} V',
$$
and the symbol $c(\cdot)$ stands for the {\it threshold of the
canonical adjunction} of the linear system $|\cdot|$,
$$
c(\Lambda)=\mathop{\rm sup}\{\varepsilon\in{\mathbb Q}_+\,\,| \,\,
D+\varepsilon K\in A^1_+(\cdot)\},
$$
$D\in\Lambda$ is an arbitrary divisor of the linear system
$\Lambda$, $K$ stands for the canonical class of the variety,
$A^1_+(\cdot)\subset A^1(\cdot)\otimes{\mathbb R}$ means the
closed cone of effective cycles on the variety under
consideration. The variety $V$ is said to be {\it birationally
superrigid}, if the inequality (\ref{0.1}) is always true for
$\chi^*=\mathop{\rm id}\nolimits_V$.

Let $\pi\colon V\to {\mathbb P}^1$ be a fibration into rationally
connected varieties, where $V$ has ${\mathbb Q}$-factorial
terminal singularities. (By the theorem of Graber-Harris-Starr
[5], in this case the variety $V$ is itself automatically
rationally connected.) The next question is of crucial importance
for understanding birational geometry of $V$: \vspace{1cm}

\parshape=1
3cm 10cm \noindent are there other (that is, different from the
original map $\pi\colon V\to {\mathbb P}^1$) structures of a
rationally connected fibration on $V$? \vspace{1cm}

Assume that $V$ is non-singular, $\mathop{\rm Pic} V={\mathbb Z}
K_V\oplus \pi^*\mathop{\rm Pic} {\mathbb P}^1$ and the following
condition holds:

\begin{equation}
\label{0.2} K_V\not\in \mathop{\rm Int} A^1_+V.
\end{equation}
Let $F_t=\pi^{-1}(t)$ be the fiber over a point $t\in{\mathbb
P}^1$, $F\in\mathop{\rm Pic} V$  the class of a fiber. The
condition (\ref{0.2}) means that if $D\sim -nK_V+lF$ is an
effective divisor on $V$, then $l\in{\mathbb Z}_+$. Conditions of
this type for three-dimensional Mori fiber spaces are discussed in
[2]. The following fact is well-known (see [17,20-22]).

{\bf Proposition 0.1.} {\it In the assumptions above let $V$ be a
birationally superrigid variety. Then:

{\rm (i)} there is only one non-trivial structure of a fibration
into rationally connected varieties on $V$, that is, the morphism
$\pi$; in other words, if $\tau\colon W\to T$ is a fibration into
rationally connected (or just uniruled) varieties and $\chi\colon
V- \,-\,\to W$ is a birational map, then $\chi$ transform fibers
into fibers, that is, the following diagram
$$
\begin{array}{rcccl}
& V & \stackrel{\chi}{- \,-\,\to} & W &  \\
\pi &  \downarrow &  & \downarrow & \tau \\
& {\mathbb P}^1 &  \stackrel{\alpha}{- \,-\,\to} & T &
\end{array}
$$
commutes for a certain map $\alpha\colon {\mathbb P}^1\to T$.

{\rm (i)} If $\tau\colon W\to T$ is another fibration of the same
type, that is, $\mathop{\rm Pic} W={\mathbb Z} K_W\oplus
\tau^*\mathop{\rm Pic} {\mathbb P}^1$, and $\chi\colon V- \,-\,\to
W$ is a (fiber-wise) birational map, then $\chi$ is an isomorphism
of fibers of general position.}

Thus the property of being birationally rigid reduces birational
geometry of the variety $V$ to biregular geometry of the fibration
$V/{\mathbb P}^1$. That is why the word ``rigidity'' has been
chosen: the variety $V$ does not admit any birational
modifications inside the natural class of Fano fibrations with the
relative Picard number one.

Nowadays quite a few classes of birationally (super)rigid Fano
varieties are known (see, for instance, [3,14,18,19,23-25]). The
known examples make it possible to conjecture that birational
(super)rigidity is a typical property in dimension 3 and higher.
Much less is known about Fano fibrations, their birational
geometry is harder to investigate. A brief history of the theory
of birational rigidity for Fano fibrations see below in Sec. 0.4.
The aim of the present paper is to prove birational superrigidity
of fibrations $V/{\mathbb P}^1$, the fibers of which are Fano
double hypersurfaces of index 1 [19].

%%%%%%%%%%%%%%%%%%%%%%%%%%%%%%%%%%%%%%%%%%%%%%%%%%%%%%%%%%%%%%%%%%%
%%%%%%%%%%%%%%%%%%%%%%%%%%%%%%%%%%%%%%%%%%%%%%%%%%%%%%%%%%%%%%%%%%%
%%%%%%%%%%%%%%%%%%%%    subsection 0.2   %%%%%%%%%%%%%%%%%%%%%%%%%%

\subsection{Varieties with a pencil of double covers}

The symbol ${\mathbb P}$ stands for the projective space ${\mathbb
P}^{M+1}$ over the field of complex numbers ${\mathbb C}$. Let
$$
{\cal G}={\mathbb P}(H^0({\mathbb P},{\cal O}_{{\mathbb P}}(m)))
$$
be the space of all Fano hypersurfaces of degree $m$, $3\leq m\leq
M-1$, ${\cal W}$ the space of all hypersurfaces of degree $2l$ in
${\mathbb P}$, where $m+l=M+1$. Let
$$
{\cal F}=\{F\,\, |\,\, \sigma\colon F \stackrel{2:1}{\to} G\}
$$
be the class of algebraic varieties realized as double covers of
hypersurfaces $G\in{\cal G}$ branched over $W\cap G$, $W\in{\cal
W}$. Set ${\cal F}_{sm}\subset {\cal F}$ to be the set of smooth
double hypersurfaces corresponding to pairs $(G,W\cap G)$ of
smooth varieties. Obviously, $F\in {\cal F}_{sm}$ is a smooth Fano
variety of index 1 with the Picard group $\mathop{\rm Pic}
F={\mathbb Z} K_F$. Let ${\cal F}^{reg}_{sm}\subset {\cal F}_{sm}$
be the smooth subset consisting of varieties $F\in {\cal F}_{sm}$
satisfying the {\it regularity condition} of Sec. 1.3,1.4 below
(which is identical to the regularity condition of Sec. 1.3 in
[19]). Recall that in [19] the following fact was proved.
\vspace{0.3cm}

{\bf Theorem A.} {\it {\rm(i)} Any variety $F\in {\cal
F}^{reg}_{sm}$ is birationally superrigid.

{\rm(ii)}  The set ${\cal F}^{reg}_{sm}$ is non-empty. Moreover,
the following estimate holds:
$$
\mathop{\rm codim}\nolimits_{{\cal F}_{sm}} ({\cal
F}_{sm}\setminus {\cal F}^{reg}_{sm})\geq 2.
$$
}

Set ${\cal F}_{sing}={\cal F}\setminus {\cal F}_{sm}$,
$\mathop{\rm codim}\nolimits_{{\cal F}}{\cal F}_{sing}=1$. Let
${\cal F}^{reg}_{sing}$ be the open subset in ${\cal F}_{sing}$,
consisting of all singular double hypersurfaces satisfying the
regularity condition of Sec. 1.3,1.4 below. We note in Sec.
1.3,1.4 that the following inequality holds:
$$
\mathop{\rm codim}\nolimits_{{\cal F}}({\cal F}_{sing}\setminus
{\cal F}^{reg}_{sing})\geq 2.
$$
In the present paper we study Fano fibrations $V/{\mathbb P}^1$,
each fiber $F_t=\pi^{-1}(t)$, $t\in{\mathbb P}^1$ of which is a
variety from the family ${\cal F}$. Set
$$
{\cal F}^{reg}={\cal F}^{reg}_{sm}\cup {\cal F}^{reg}_{sing}.
$$
By what was said above, $\mathop{\rm codim}\nolimits_{{\cal
F}}({\cal F} \setminus {\cal F}^{reg})\geq 2$. Since the fibration
$V/{\mathbb P}^1$ can be looked at as a morphism ${\mathbb P}^1
\to {\cal F}$, that is, a curve in ${\cal F}$, for a general
variety $V/{\mathbb P}^1$ we get:
$$
F_t\in {\cal F}^{reg}
$$
for all points $t\in {\mathbb P}^1$. If this is the case, we say
that the fibration $V/{\mathbb P}^1$ is {\it regular}. A general
construction of regular Fano fibrations $V/{\mathbb P}^1$ is
described below in Sec. 1.2.

%%%%%%%%%%%%%%%%%%%%%%%%%%%%%%%%%%%%%%%%%%%%%%%%%%%%%%%%%%%%%%%%%%%
%%%%%%%%%%%%%%%%%%%%%%%%%%%%%%%%%%%%%%%%%%%%%%%%%%%%%%%%%%%%%%%%%%%
%%%%%%%%%%%%%%%%%%%%    subsection 0.3   %%%%%%%%%%%%%%%%%%%%%%%%%%
\subsection{The main result}

{\bf Theorem 1.} {\it Assume that a regular fibration $V/{\mathbb
P}^1$ satisfies the $K^2$-condition:
$$
K^2_V\not\in \mathop{\rm Int} A^2_+ V,
$$
where $A^2_+ V\subset A^2 V\otimes{\mathbb R}$ is the closed cone
of effective cycles of codimension two. The the fibration
$V/{\mathbb P}^1$ is birationally superrigid.}

The symbol $A^i \sharp$ stands, as usual, for the group of classes
of codimension $i$ cycles on the variety $\sharp$ modulo numerical
equivalence.

{\bf Corollary 0.1.} {\it {\rm (i)} For a fibration $V/{\mathbb
P}^1$ of general position the following equality holds:
$$
\mathop{\rm Bir} V=\mathop{\rm Aut} V= {\mathbb Z}/2{\mathbb
Z}=\{\mathop{\rm id},\tau\},
$$
where $\tau\in \mathop{\rm Aut} V$ is the Galois involution of the
double cover $V/Q$.

{\rm(ii)} The variety $V$ is non-rational.}

Part (i) follows from Proposition 0.1 and Theorem A. Part (ii) is
obvious.

We prove Theorem 1 in a few steps. Below in Sec. 1.1 we formulate
a sufficient condition of birational superrigidity for an
arbitrary Fano fibration $V/{\mathbb P}^1$  in terms of numerical
geometry of fibers (Theorem 2). Essentially this fact was proved
in [17,20], although the papers that we have just mentioned
discussed Fano fibrations of a certain particular type. We will
not repeat these arguments here, just making reference to [17,20].

Now to prove Theorem 1 we need to check that the fibers of the
fibration $V/{\mathbb P}^1$, that is, the regular Fano
hypersurfaces of index 1 (in the sense of the regularity
conditions formulated below in Sec. 1.3,1.4) satisfy the
conditions of Theorem 2. This verification makes our proof. It is
carried out in Sections 2 and 3. Birational geometry of varieties
with a pencil of Fano double covers that do not satisfy the
$K^2$-condition will be studied in the next paper, the second part
of the present research.

{\bf Remark.} Superrigidity of Fano fibrations  $V/{\mathbb P}^1$,
the fibers of which are double spaces ($m=1$) and double quadrics
($m=2$) of index 1, is proved in [17] and [21,22], respectively,
and for this reason these varieties are not considered in this
paper.

%%%%%%%%%%%%%%%%%%%%%%%%%%%%%%%%%%%%%%%%%%%%%%%%%%%%%%%%%%%%%%%%%%%
%%%%%%%%%%%%%%%%%%%%%%%%%%%%%%%%%%%%%%%%%%%%%%%%%%%%%%%%%%%%%%%%%%%
%%%%%%%%%%%%%%%%%%%%    subsection 0.4   %%%%%%%%%%%%%%%%%%%%%%%%%%
\subsection{Historical remarks}

Investigating structures of a fibration into rationally connected
(uniruled) varieties is a very old subject. The classical proof of
the Noether theorem on the Cremona group of the plane, presented
by Yu.I.Manin in [1], can be looked at in this way: step by step a
certain pencil of rational curves on ${\mathbb P}^2$ is modified,
thus one goes over from one structure of a ${\mathbb
P}^1$-fibration on the plane to another. Birational geometry of
varieties of dimension higher than two, on which there are a lot
of various structures of a rationally connected fibration, is very
hard to study. Apart from a few exceptional types, these varieties
are still out of reach for the modern technique. However, on
rationally connected varieties which are in a certain sense
general there is only one structure of a rationally connected
fibration with the minimality condition that it is equivalent to a
Fano fibration with the relative Picard number one. This is the
very phenomenon of birational rigidity.

In the modern birational geometry birationally rigid varieties
first come to the light in the papers of Yu.I.Manin in the form of
del Pezzo surfaces over non-closed fields [15,16]. The first
theorems on birational rigidity of non-trivial rationally
connected fibrations were proved by V.A.Iskovskikh as theorems on
uniqueness of a pencil of rational curves for some surfaces over
non-closed fields. These theorems continue the above-mentioned
work of Yu.I.Manin. The study of the absolute and relative cases
in dimension two over a non-closed field prepared the basis for
working in higher dimensions.

In the classical paper of V.A.Iskovskikh and Yu.I.Manin [14] the
{\it test class} technique was developed that made it possible to
prove (in the modern terminology) birational superrigidity of the
smooth three-dimensional quartic $V_4\subset{\mathbb P}^4$
(actually, in [14] birational superrigidity of the double space
branched over a sextic was proved as well and also the crucial
step was made in the proof of birational rigidity of the double
quadric of index one [11]). After that attempts were made to use
this technique to obtain similar results in the relative case for
fibrations over a non-trivial base. For one class of varieties the
study was successful: V.G.Sarkisov's theorem proves that the given
structure of a conic bundle is unique provided the discriminant
divisor is sufficiently big [26,27]. The proof of Sarkisov's
theorem is based on the following two technical principles:

(1) the test class technique of V.A.Iskovskikh and Yu.I.Manin,

(2) the fiber-wise modifications.

The possibility of making fiber-wise modifications of conic
bundles remaining at the same time in the class of smooth
varieties is an exclusive property of these varieties. In a sense
this special feature comes from the fact that the group of
automorphisms of the fiber $\mathop{\rm Aut}{\mathbb P}^1$ is very
big: say, for a typical Fano variety of dimension higher than two
this group is finite. Thus there is no hope to use similar
arguments for higher dimensional Fano fibrations.

After Sarkisov's papers [26,27] had been published, there remained
only one class of rationally connected three-folds, birational
geometry of which was a terra incognita, that is, the class of
fibrations into del Pezzo surfaces over ${\mathbb P}^1$. The
attempts to use fiber-wise modifications similar to Sarkisov's
theorem proved unsuccessful, since immediately converted the
variety under consideration into a singular one and, moreover, the
acquired singularities were out of control. However, the
above-mentioned test class technique also refused to work. Since
mid-80s and up to mid-90s attempts were made to construct at least
some examples of three-dimensional del Pezzo fibrations similar to
Sarkisov's rigid conic bundles, but without any success. This
activity was summed up in [12]: the only outcome of the almost
decade-long work and an immense amount of completed computations
were some conjectures --- and no essential progress in their
proof. One can see from [12] that there was no understanding why
the test class technique that works so impeccably in the absolute
case (the three-dimensional quartic [14]) does not allow a single
step forward in the case of del Pezzo fibrations: the test class
simply refused to be constructed.

The situation changed radically when the paper [17] appeared. It
became immediately clear why the test class technique refused to
generalize to the relative case: as it turned out, the desired
class just did not exist. For the three-dimensional quartic the
test class technique is equivalent to the technique of counting
multiplicities introduced in [17,18,23]. However, the technique of
counting multiplicities is much more flexible, since it describes
properties of a certain effective cycle of codimension two (the
self-intersection of the moving linear system defining the
birational map under consideration), whereas the test class gives
just a number, the intersection number of this cycle with the test
class. In the relative case, when the base of the fibration is
non-trivial, any effective cycle can be decomposed into the
vertical and horizontal components. Informally speaking, each of
them requires its own test class.

The methods developed in the paper [17] were later used for
proving birational rigidity of big classes of higher-dimensional
Fano fibrations [20-22]. In these papers (and in the present paper
as well) birational rigidity is derived from the $K^2$-condition.
However when the $K^2$-condition is somewhat weakened the methods
of these papers still work well and make it possible to give a
complete description of birational geometry of the variety under
consideration. See [29,30] and the series of papers [6-8], where
birational rigidity is proved for a few classes of del Pezzo
fibrations over ${\mathbb P}^1$. These classes were not considered
in [17] because they do not satisfy the $K^2$-condition. However,
when the deviation from the $K^2$-condition grows too strong, the
methods fail to work.

Note also that in spite of the progress in the general theory of
factorization of birational maps between three-fold Mori fiber
spaces (the Sarkisov program [4,28]), all attempts either to
improve Sarkisov's results or to prove the rationality criterion
for conic bundles have been unsuccessful up to this day  [13]. See
the recent paper [2] on this point. We will discuss it in the next
papers.

%%%%%%%%%%%%%%%%%%%%%%%%%%%%%%%%%%%%%%%%%%%%%%%%%%%%%%%%%%%%%%%%%%%
%%%%%%%%%%%%%%%%%%%%%%%%%%%%%%%%%%%%%%%%%%%%%%%%%%%%%%%%%%%%%%%%%%%
%%%%%%%%%%%%%%%%%%%%    subsection 0.5   %%%%%%%%%%%%%%%%%%%%%%%%%%
\subsection{Acknowledgements}

One part of the present research (corresponding to the second
section of this paper) was carried out by the author during his
stay at the University of Bayreuth in 2001 as a Humboldt Research
Fellow. The crucial step (investigation of singularities coming
from a double point on the branch divisor of a fiber, making the
contents of the third section of the paper) was made during my
work at Max-Planck-Institut f\" ur Mathematik in Bonn in 2003. The
author is very grateful to Alexander von Humboldt Stiftung,
Mathematisches Institut der Universit\" at Bayreuth (in the first
place, to Prof. Th.Peternell) and Max-Planck-Institut f\" ur
Mathematik in Bonn for hospitality, the excellent conditions of
work and general support.

%%%%%%%%%%%%%%%%%%%%%%%%%%%%%%%%%%%%%%%%%%%%%%%%%%%%%%%%%%%%%%%%%%%
%%%%%%%%%%%%%%%%%%%%%%%%%%%%%%%%%%%%%%%%%%%%%%%%%%%%%%%%%%%%%%%%%%%
%%%%%%%%%%%%%%%%%%%%%%%%%%%%%%%%%%%%%%%%%%%%%%%%%%%%%%%%%%%%%%%%%%%
%%%%%%%%%%%%%%%%%%%%%%%%%%%%%%%%%%%%%%%%%%%%%%%%%%%%%%%%%%%%%%%%%%%
%%%%%%%%%%%%%%%%%%%    section 1                     %%%%%%%%%%%%%%
\section{The method of maximal singularities and the regularity
conditions}

\subsection{A criterion of birational rigidity}

Let $\pi\colon V\to{\mathbb P}^1$ be a smooth standard Fano
fibration, that is, $V$ be a smooth variety with
$$
\mathop{\rm Pic} V={\mathbb Z} K_V\oplus {\mathbb Z} F,
$$
where $F$ is the class of a fiber. Define the {\it degree} of a
horizontal subvariety $Y\subset V$, $\pi(Y)={\mathbb P}^1$, by the
formula
$$
\mathop{\rm deg}\nolimits Y= (Y\cdot F\cdot (-K_V)^{\mathop{\rm
dim} Y-1}),
$$
and the degree of a vertical subvariety $Y\subset \pi^{-1}(t)$ by
the formula
$$
\mathop{\rm deg}\nolimits Y=(Y\cdot (-K_V)^{\mathop{\rm dim} Y}).
$$
By this definition the degree of the variety $V$ itself coincides
with the degree of a fiber, $\mathop{\rm deg}\nolimits
V=\mathop{\rm deg}\nolimits F$.

Smooth Fano fibrations, the fibers of which are complete
intersections in weighted projective spaces, satisfy also the
following property: their fibers have at most isolated
singularities.

We say that the Fano fibration $V/{\mathbb P}^1$ satisfies
\vspace{0.5cm}

\parshape=1
1cm 12cm \noindent {\it condition} (v), if for any irreducible
vertical subvariety $Y$ of codimension 2, $Y\subset
\pi^{-1}(t)=F_t$, and any smooth point $o\in F_t$ the following
estimate holds:
$$
\frac{\mathop{\rm mult}\nolimits_o}{\mathop{\rm deg}} Y \leq
\frac{2}{\mathop{\rm deg}\nolimits V};
$$
\vspace{0.5cm}

\parshape=1
1cm 12cm \noindent {\it condition} (vs), if for any vertical
subvariety $Y\subset F_t$ of codimension 2 (with respect to $V$,
that is, a prime divisor on $F_t$), any singular point $o\in F_t$
and any infinitely near point $x\in \widetilde F_t$, where
$\varphi\colon \widetilde F_t\to F_t$ is the blow up of the point
$o$, $\varphi(x)=o$, $\widetilde Y\subset \widetilde F_t$ the
strict transform of the subvariety $Y$ on $\widetilde F_t$, the
following estimates hold:
$$
\frac{\mathop{\rm mult}\nolimits_o}{\mathop{\rm deg}} Y \leq
\frac{4}{\mathop{\rm deg}\nolimits V}, \quad \frac{\mathop{\rm
mult}\nolimits_x \widetilde Y }{\mathop{\rm deg}\nolimits Y}  \leq
\frac{2}{\mathop{\rm deg}\nolimits V};
$$
\vspace{0.5cm}

\parshape=1
1cm 12cm \noindent {\it condition} (h), if for any horizontal
subvariety $Y$ of codimension 2 and any point $o\in Y$ the
following estimate holds
$$
\frac{\mathop{\rm mult}\nolimits_o}{\mathop{\rm deg}} Y \leq
\frac{4}{\mathop{\rm deg}\nolimits V}.
$$

Assume that $\mathop{\rm dim} V\geq 4$ and the variety $V$
satisfies the condition
$$
A^2 V={\mathbb Z} K^2_V\oplus {\mathbb Z} H_F,
$$
where $H_F=(-K_V\cdot F)$ and a fiber $F=F_t\subset V$ of general
position satisfies the condition $A^2 F={\mathbb Z} (H_F\cdot
H_F)_F$. Set $A^2_{\mathbb R} V=A^2 V\otimes {\mathbb R} \cong
{\mathbb R}^2$ and define the cone of effective cycles $A^2_+
V\subset  A^2_{\mathbb R} V$ as the closure (in the real topology)
of the set
$$
\{\lambda\Delta\,\, |\,\, \lambda\in{\mathbb R}_+,
\Delta\,\,\mbox{is the class of an effective cycle}\}.
$$

{\bf Definition 1.1} We say that the Fano fibration $V/{\mathbb
P}^1$ satisfies the $K^2$-{\it condition}, if
$$
K^2_V\not\in \mathop{\rm Int} A^2_+ V.
$$

{\bf Remark.} It is easy to see that the $K^2$-condition is
equivalent to the following claim: for any $a\geq 1$ and $b\geq 1$
the class
$$
\Delta(a,b)=aK^2_V-b H_F
$$
is not effective. Indeed, $H_F\in A^2_+ V$, so that $K^2_V\in
\mathop{\rm Int} A^2_+ V$ if and only if $\Delta(N,1) \in A^2_+ V$
for some $N\geq 1$. This implies immediately that both conditions
are equivalent.

{\bf Theorem 2.} {\it Assume that the smooth standard Fano
fibration $V/{\mathbb P}^1$ satisfies the $K^2$-condition and the
conditions {\rm (v)}, {\rm (vs)} and {\rm (h)}. Then $V/{\mathbb
P}^1$ is birationally superrigid.}

For the {\bf proof} see [17,20].

%%%%%%%%%%%%%%%%%%%%%%%%%%%%%%%%%%%%%%%%%%%%%%%%%%%%%%%%%%%%%%%%%%%
%%%%%%%%%%%%%%%%%%%%%%%%%%%%%%%%%%%%%%%%%%%%%%%%%%%%%%%%%%%%%%%%%%%
%%%%%%%%%%%%%%%%%%%    subsection 1.2    %%%%%%%%%%%%%%%%%%%%%%%%%%
\subsection{An explicit construction of the fibration
$V/{\mathbb P}^1$}

Let us describe an explicit construction of regular fibrations
$V/{\mathbb P}^1$. For each fiber $F\in{\cal F}^{reg}$ (singular
or smooth) the anticanonical linear system $|-K_F|$ determines
precisely the double cover $\sigma_F\colon F\to G\subset{\mathbb
P}$. For this reason, $\pi_*{\cal O}(-K_V)$ is a locally free
sheaf of rank $M+2$ on ${\mathbb P}^1$. It gives a locally trivial
${\mathbb P}$-fibration over ${\mathbb P}^1$. The variety $V$ is
realized as a double cover of a smooth divisor $Q$ on ${\mathbb
P}(\pi_*{\cal O}(-K_V))$. Namely, let
$$
{\cal E}=\mathop{\bigoplus}\limits^{M+1}_{i=0} {\cal O}_{{\mathbb
P}^1}(a_i)
$$
be a locally free sheaf, normalized by the condition that
$$
a_0=0\leq a_1\leq\dots\leq a_i\leq a_{i+1}\leq\dots\leq a_{M+1}.
$$
In particular, ${\cal E}$ is generated by global sections. Set
$X={\mathbb P}({\cal E})$ to be its ${\bf Proj}$ in the sense of
Grothendieck, $\pi_X\colon X\to{\mathbb P}^1$ the natural
projection, ${\cal L}_X$ the tautological sheaf, $Q\subset X$ a
smooth divisor on $X$, corresponding to a section
$$
s_Q\in H^0\left(X,{\cal L}^{\otimes m}_X\otimes \pi^*_X{\cal
O}_{{\mathbb P}^1} (a_Q) \right),
$$
$a_Q\in{\mathbb Z}_+$. The symbol $\pi_Q\colon Q\to{\mathbb P}^1$
stands for the projection $\pi_X|_Q$. Obviously, $Q/{\mathbb P}^1$
is a smooth fibration into Fano hypersurfaces of degree $m$ in
${\mathbb P}$. Let $W\subset X$ be an irreducible hypersurface,
corresponding to a section
$$
s_W\in H^0\left(X,{\cal L}^{\otimes 2l}_X\otimes \pi^*_X{\cal
O}_{{\mathbb P}^1} (2a_W) \right),
$$
$a_W\in{\mathbb Z}_+$, whereas $W_Q=W\cap Q$ is a smooth divisor
on $Q$. We denote the fiber $\pi^{-1}_Q(t)$ over a point
$t\in{\mathbb P}^1$  by the symbol $G_t$ (or just $G$, when it is
clear which point is meant or when it is inessential). Finally,
set
$$
\sigma\colon V\to Q
$$
to be the double cover, branched over $W_Q$. The natural
projection onto ${\mathbb P}^1$ will be denoted by $\pi$, the
fiber $\pi^{-1}(t)$ by the symbol $F_t$ (or just $F$). It is easy
to see that
$$
\mathop{\rm Pic} V={\mathbb Z} K_V \oplus {\mathbb Z} F
$$
and up to twisting by an invertible sheaf ${\cal O}_{{\mathbb
P}^1}(k)$, $k\in{\mathbb Z}$, the sheaves ${\cal E}$ and
$\pi_*{\cal O}(-K_V)$ on ${\mathbb P}^1$ coincide.

More precisely, let $L_X\in \mathop{\rm Pic} X$ be the class of
the tautological sheaf ${\cal L}_X$, $L_Q=L_X|_Q$ its restriction
to $Q$, so that
$$
\mathop{\rm Pic} Q={\mathbb Z} L_Q \oplus {\mathbb Z} G.
$$
Set $L_V=\sigma^* L_Q$. It is easy to see that
$$
K_V=-L_V+(a_1+\dots+a_M-2+a_Q+a_M)F.
$$
By the Lefschetz theorem
$$
A^2 V={\mathbb Z} K^2_V\oplus {\mathbb Z} H_F,
$$
where $H_F=(-K_V\cdot F)$ is the class of a hyperplane section.
The symbol $H_F$ is used in the present paper in two different
meanings: as a class of codimension two on $V$ and as the
hyperplane section of the fiber, that is, an element of $A^1F$.
Every time it is clear which of the two concepts is meant.

It is easy to compute that
$$
(K^2_V\cdot L^{M-1})=2m (4-a_1-\dots -a_{M+1}-a_Q-a_W)+2a_Q.
$$
Since $(H_F\cdot L^{M-1})=2m$ and the linear system $|L_V|$ is
free, the inequality $(K^2_V\cdot L^{M-1})\leq 0$ implies, that
$K^2_V\not\in \mathop{\rm Int} A^2_+ V$, where $A^2_+ V\subset A^2
V\otimes {\mathbb R}$ is the closed cone of effective cycles of
codimension two.

%%%%%%%%%%%%%%%%%%%%%%%%%%%%%%%%%%%%%%%%%%%%%%%%%%%%%%%%%%%%%%%%%%%
%%%%%%%%%%%%%%%%%%%%%%%%%%%%%%%%%%%%%%%%%%%%%%%%%%%%%%%%%%%%%%%%%%%
%%%%%%%%%%%%%%%%%%%    subsection 1.3   %%%%%%%%%%%%%%%%%%%%%%%%%%%
\subsection{The regularity conditions outside the branch divisor}

Let $\sigma\colon F\to G\subset {\mathbb P}$ be a Fano double
hypersurface of index 1, $F\in{\cal F}$. The variety $F$ is
realized as a complete intersection of codimension two in the
weighted projective space
$$
{\mathbb P}(\underbrace{1,1,\dots,1,}_{M+2} l),
$$
see [19]: $F$ is of type $m\cdot 2l$ and given by the pair of
equations
$$
\tilde f(x_0,\dots,x_{M+1})=0,\quad u^2=\tilde
g(x_0,\dots,x_{M+1}),
$$
where $x_*$ are the coordinates of weight 1, $u$ is the coordinate
of weight $l$, $\tilde f$ is the equation of the hypersurface
$G\subset {\mathbb P}={\mathbb P}(1,\dots,1)$, $\tilde g$ is the
equation of the hypersurface $W\cap {\mathbb P}$.

Let $o\in F$ be an arbitrary point. First of all, we draw the
reader's attention to the following obvious fact:
$$
o\neq   (\underbrace{0,0,\dots,0,}_{M+2} 1).
$$
Thus we may assume that the point $o$ lies in one of the standard
affine charts ${\mathbb A}^{M+2}$ with the coordinates
$$
z_i=x_i/x_0,\quad i=1,\dots,M+1,\quad y=u/x_0^l
$$
and its $z_*$-coordinates are $(0,\dots,0)$. With respect to the
coordinate system $(z_*,y)$ the affine part of the variety $F$ is
given by the pair of equations
$$
f=q_1+\dots+q_M=0,\quad y^2=g=w_0+\dots+w_{2l},
$$
where $q_i$ and $w_j$ are homogeneous polynomials in $z_*$ of
degrees $i$ and $j$, respectively. Set $p=\sigma(o)\in G$. The
point $p$ lies on the branch divisor $W$ if and only if $w_0=0$.
If $p\not\in W$, then we normalize the second equation and assume
that $w_0=1$.

Let us formulate first the regularity conditions outside the
branch divisor. In this case the fiber $F$ is given with respect
to the affine coordinate system $(z_*,y)$ with the origin of the
$z_*$-system at $p=\sigma(o)$ by the equations
$$
\left\{
\begin{array}{l}
f=q_a+\dots+q_m=0, \\
y^2=g=1+w_1+\dots+w_{2l},
\end{array}
\right.
$$
where $a\geq 1$.

Set
$$
\sqrt{g}=(1+w_1+\dots+w_{2l})^{1/2}=
1+\sum^{\infty}_{i=1}\gamma_i(w_1+\dots+w_{2l})^i=
$$
$$
=1+ \sum^{\infty}_{i=1}\Phi_i(w_1,\dots,w_{2l}),
$$
where $\Phi_i(w_1(z_*),\dots,w_{2l}(z_*))$ are homogeneous in
$z_*$ of degree $i\geq 1$,
$$
\gamma_i=(-1)^{i-1}\frac{(2i-3)!!}{2^i i!}=
(-1)^{i-1}\frac{(2i-3)!}{2^{2i-2} i! (i-2)!}
$$
is the standard $i$-th coefficient of the Taylor expansion of the
function $(1+s)^{1/2}$ at the point $s=0$. Obviously,
$$
\Phi_i(w_*)=w_i+A_i(w_1,\dots,w_{i-1})
$$
for $i\leq 2l$. For $i\geq 1$ set
$$
[\sqrt{g}]_i=1+\sum^i_{j=1}\Phi_i(w_*),\quad
g^{(i)}=g-[\sqrt{g}]^2_i.
$$
It is easy to see that the first non-zero component of the
polynomial $g^{(i)}$ is of degree $i+1$. More precisely, this
component is equal to
$$
g_{i+1}=2\Phi_{i+1}(w_1(z_*),\dots,w_{i+1}(z_*)).
$$
\vspace{0.5cm}

{\bf The regularity condition at a smooth point} $p\in G$
(R1.1):\vspace{0.5cm}

The sequence
$$
q_1,\dots,q_m,g_{l+1},\dots,g_{2l-1}
$$
is regular in ${\cal O}_{p,{\mathbb P}}$. Here
$a=1$.\vspace{0.5cm}

{\bf  The regularity condition at a double point} $p\in G$
(R1.2):\vspace{0.5cm}

If $2l\geq m+1$, then the system of $M-1$ homogeneous polynomials
$$
q_2,\dots,q_m,g_{l+1},\dots,g_{2l-1},
$$
whereas if $2l\leq m$, then the system of homogeneous polynomials
$$
q_2,\dots,q_{m-1},g_{l+1},\dots,g_{2l}
$$
defines a curve in ${\mathbb P}^M={\mathbb P}(T_p {\mathbb P})$,
neither component of which is contained in a hyperplane.

Furthermore, the system of $M$ homogeneous equations
\begin{equation}
\label{1.1} q_2=\dots=q_{m-1}=g_{l+1}=\dots=g_{2l}=0
\end{equation}
defines a non-zero subscheme $Z_*$ in ${\mathbb P}^M$, such that
for any hyperplane $P\subset{\mathbb P}^M$
$$
\mathop{\rm deg}\nolimits (P\cap
Z_*)<\lambda_{m,l}=\frac{m!(2l-1)! }{6 (l-1)!}
$$
for $m\geq 4$ and
$$
\mathop{\rm deg}\nolimits (P\cap
Z_*)<\lambda_{3,l}=12\frac{(2l-1)! }{(l+1)!} (l-2)
$$
for $m=3$. If the scheme $Z_*$ is reduced, then this condition
means simply that any set of $\lambda_{m,l}$ points is not
contained in a hyperplane.\vspace{0.5cm}

{\bf Remark.} Since $w_0=y(0)=1$, in a neighborhood of the
singular point $o\in F$ the equations
$$
y-[\sqrt{g}]_i=0\quad\mbox{and}\quad \sigma^* g^{(i)}=0
$$
define the same divisor. Consider the system of equations
(\ref{1.1}) on the fiber $F$ (and not on the projectivized tangent
space ${\mathbb P}(T_p {\mathbb P})$). The system defines an
effective 1-cycle $C_*$ on $F$. By construction, its degree is
equal to
$$
\mathop{\rm deg}\nolimits C_*=2m!\frac{(2l-1)! }{(l-1)!},
$$
whereas its multiplicity at the point $o\in F$ satisfies the
estimate
\begin{equation}
\label{1.2} \mathop{\rm mult}\nolimits_o C_*\geq m! \frac{(2l)!
}{l!}=\mathop{\rm deg}\nolimits C_*,
\end{equation}
so that what we actually have in (\ref{1.2}) is an equality and
$C_*$ is an algebraic sum of lines on $F$, that is, curves of the
form $L\ni o$, the image $\sigma(L)\subset {\mathbb P}$ of which
is a line, and moreover the morphism $\sigma\colon L\to\sigma(L)$
is an isomorphism. Considering the zero-dimensional scheme $Z_*$
as an effective zero-dimensional cycle, we get by construction:
$$
Z_*={\mathbb P}(T_oC_*).
$$
In particular, for any hyperplane $P\subset{\mathbb P}$ the
one-dimensional part of the scheme
$$
\{q_2=\dots=q_{m}=g_{l+1}=\dots=g_{2l}=0\}\cap \sigma^{-1}(P)
$$
is of degree not higher than $\lambda_{m,l}-1$. In other words, if
all components of the cycle $C_*$ are of multiplicity 1, then no
more than $\lambda_{m,l}-1$ of these lines are contained in
$\sigma^{-1}(P)$.

%%%%%%%%%%%%%%%%%%%%%%%%%%%%%%%%%%%%%%%%%%%%%%%%%%%%%%%%%%%%%%%%%%%
%%%%%%%%%%%%%%%%%%%%%%%%%%%%%%%%%%%%%%%%%%%%%%%%%%%%%%%%%%%%%%%%%%%
%%%%%%%%%%%%%%%%%%%    subsection 1.4   %%%%%%%%%%%%%%%%%%%%%%%%%%%
\subsection{The regularity conditions on the branch divisor}

In this case the variety  $F$ is given with respect to the affine
coordinate system $(z_*,y)$ by the system of equations
$$
\left\{
\begin{array}{l}
f=q_1+\dots+q_m=0, \\
y^2=g=w_1+\dots+w_{2l}.
\end{array}
\right.
$$
\vspace{0.5cm}

{\bf The regularity condition at a smooth point} $o\in F$
(R2.1):\vspace{0.5cm}

\noindent the sequence of homogeneous polynomials
$$
q_1,\dots,q_m
$$
is regular in ${\cal O}_{p,{\mathbb P}}$ and the quadratic form
$q_2$ does not vanish identically on the plane $\{q_1=w_1=0\}$.

Note that since the point $o\in F$ is smooth, this plane is of
codimension exactly two, that is, the linear forms $q_1$ and $w_1$
are linearly independent: the plane $\{q_1=w_1=0\}$ is the tangent
plane to the branch divisor $W\cap G$ of the morphism
$\sigma_F$.\vspace{0.5cm}

{\bf The regularity condition at a double point} $o\in F$
(R2.2):\vspace{0.5cm}

\noindent In this case we have the double cover $\sigma_F\colon
F\to G$, branched over the divisor $W_G=W\cap G$. The first
regularity condition is smoothness of the hypersurface $G$ at the
point $p=\sigma(o)$, that is, $q_1\neq 0$. Furthermore, the
divisor $W_G$ should have at the point $p$ a non-degenerate
quadratic singularity:
$$
w_1=\lambda q_1,
$$
$\lambda\in{\mathbb C}$. For convenience of notations assume that
$q_1=z_{M+1}$. The quadratic polynomial
$$
{\bar w}_2=w_2|_{\{z_{M+1}=0\}}
$$
is of the maximal rank. Let $E_G\cong {\mathbb P}^{M-1}$ be the
exceptional divisor of the blow up $\varphi_G\colon \widetilde G
\to G$ of the point $p$. Take $z_1,\dots,z_M$ for homogeneous
coordinates on $E_G$ and set
$$
W_E=\{ {\bar w}_2=0\}.
$$
It is a non-singular quadratic hypersurface in $E_G$. Denote by
the symbol ${\bar q}_i$ the restriction of the homogeneous
polynomial $q_i$ onto the hyperplane $q_{M+1}=0$. Now the
remaining part of the condition (R2.2) looks as
follows:\vspace{0.5cm}

\noindent the system of homogeneous equations
$$
{\bar q}_2=\dots={\bar q}_m=0
$$
defines in $E_G\cong {\mathbb P}^{M-1}_{(z_1:\dots :z_M)}$ an
irreducible subvariety $Z_{2\cdot \dots \cdot m}$, which is an
irreducible reduced complete intersection of codimension $(m-1)$.
The quadric
$$
{\bar q}_2=0
$$
is smooth and distinct from $W_E$.\vspace{0.5cm}

{\bf Definition 1.2.} A Fano double hypersurface $F\in {\cal F}$
is {\it regular}, if each smooth point on it is regular in the
sense of the corresponding condition (R1.1) or (R2.1) and each of
its singular points is regular in the sense of the corresponding
condition (R1.2) or (R2.2). Notation: $F\in{\cal F}^{reg}$.

The conditions (R1.1) and (R2.1) coincide with the regularity
conditions of the paper [19] (Definitions 1 and 2 in Sec. 1.3). In
[19, Sec. 4.3] it was shown that non-regular smooth double spaces
form a closed subset of codimension at least two in the set of all
smooth double hypersurfaces ${\cal F}_{sm}$. Moreover, it follows
from the computations of Sec. 4.3 in [19] that the set of Fano
double hypersurfaces $F$ with at least one {\it smooth}
non-regular point $o\in F$ is of codimension at least two in
${\cal F}$. Thus a general singular double hypersurface $F\in{\cal
F}_{sing}$ has exactly one singular point whereas all its smooth
points are regular. The singular point $o\in F$ is a
non-degenerate double point. If $p=\sigma(o)\not\in W_G$, then the
fact that the condition (R1.2) is open implies that in a
neighborhood of $F\in{\cal F}$ the following estimate holds
\begin{equation}
\label{1.3} \mathop{\rm codim}\nolimits_{{\cal F}_{sing}}({\cal
F}_{sing} \setminus {\cal F}^{reg}_{sing})\geq 1
\end{equation}
and thus
\begin{equation}
\label{1.4} \mathop{\rm codim}\nolimits_{{\cal F}}({\cal F}_{sing}
\setminus {\cal F}^{reg}_{sing})\geq 2.
\end{equation}
If $p=\sigma(o)\in W_G$, then in a similar way the fact that the
condition (R2.2) is open implies the estimate (\ref{1.3}) in a
neighborhood of $F\in{\cal F}$. Thus the estimates (\ref{1.3}) and
(\ref{1.4}) are global.

%%%%%%%%%%%%%%%%%%%%%%%%%%%%%%%%%%%%%%%%%%%%%%%%%%%%%%%%%%%%%%%%%%%
%%%%%%%%%%%%%%%%%%%%%%%%%%%%%%%%%%%%%%%%%%%%%%%%%%%%%%%%%%%%%%%%%%%
%%%%%%%%%%%%%%%%%%%    subsection 1.5   %%%%%%%%%%%%%%%%%%%%%%%%%%%
\subsection{Start of the proof of Theorem 1}

Let us check that the regular fibration $V/{\mathbb P}^1$
satisfies the conditions (v) and (h). Assume that the opposite
inequality holds:
$$
\frac{\mathop{\rm mult}\nolimits_o}{\mathop{\rm deg}} Y
>\frac{2}{\mathop{\rm deg}\nolimits V},
$$
where $o\in F=F_t$ is a smooth point, $Y\subset F$ is a prime
divisor. Let
$$
T=\sigma^{-1}(T_pG\cap G)
$$
be the tangent divisor, $p=\sigma(o)$. By the regularity
conditions, $\mathop{\rm mult}\nolimits_o T=2$. Since $T\subset F$
is a hyperplane section, we get $\mathop{\rm deg}\nolimits
T=\mathop{\rm deg}\nolimits V$, so that
$$
\frac{\mathop{\rm mult}\nolimits_o}{\mathop{\rm deg}} T
=\frac{2}{\mathop{\rm deg}\nolimits V}
$$
and thus $Y\neq T$. Both subvarieties $Y$, $T$ are irreducible, so
that the intersection $Y\cap T$ is of codimension two with respect
to $F$ and the effective cycle $Z=(Y\circ T)$ of the
scheme-theoretic intersection  of $Y$ and $T$ is well defined.
Obviously, the cycle $Z$ satisfies the inequality
\begin{equation}
\label{1.5} \frac{\mathop{\rm mult}\nolimits_o}{\mathop{\rm deg}}
Z
> \frac{4}{\mathop{\rm deg}\nolimits V}.
\end{equation}
However it was proved in [19], Sec. 3, that for a regular point
$o\in F$ it is impossible. This proves the condition (v).

Let us prove that the condition (h) holds. To begin with, let us
consider first the smooth case, where $o\in F$ is a smooth point.
Assume that an irreducible horizontal subvariety $Y\subset V$ of
codimension two satisfies the inequality
$$
\frac{\mathop{\rm mult}\nolimits_o}{\mathop{\rm deg}} Y
> \frac{4}{\mathop{\rm deg}\nolimits V}.
$$
Since $\pi(Y)={\mathbb P}^1$, we get $Y\neq F$, so that $Z=(Y\circ
F)$ is an effective cycle of codimension two on the fiber $F$,
satisfying the inequality (\ref{1.5}). As it was pointed out
above, this is impossible. The condition (h) is proved in the
smooth case.

Now let $o\in F$ be a double point. Arguing in the same way as in
the smooth case, let us construct the effective cycle $Z=(Y\circ
F)$ of codimension two on the fiber $F$. Since
$\displaystyle\frac{\mathop{\rm mult}\nolimits_o}{\mathop{\rm
deg}} F=2$, the cycle $Z$ satisfies the inequality
$$
\frac{\mathop{\rm mult}\nolimits_o}{\mathop{\rm deg}} Z
> \frac{8}{\mathop{\rm deg}\nolimits V}.
$$
Let us show that this is impossible. Without loss of generality
assume that $Z\subset F$ is an irreducible subvariety of
codimension two. Its image on $G$ satisfies the estimate
$$
\frac{\mathop{\rm mult}\nolimits_p}{\mathop{\rm deg}}
\sigma(Z)>\frac{4}{\mathop{\rm deg}\nolimits G}.
$$
Now if $p\in G$ is a smooth point, then the arguments of the paper
[20] (they work without any modifications for an arbitrary degree
$\mathop{\rm deg}\nolimits G\leq \mathop{\rm dim} G+1$) show that
this is impossible. If $p\in G$ is a double point, then by the
condition (R1.2) the homogeneous polynomials $q_2,\dots,q_m$ make
a regular sequence, so that the standard arguments of [20] give a
contradiction once again (see Sec. 3.1 in [20]).

This completes the proof of the condition (h).

%%%%%%%%%%%%%%%%%%%%%%%%%%%%%%%%%%%%%%%%%%%%%%%%%%%%%%%%%%%%%%%%%%%
%%%%%%%%%%%%%%%%%%%%%%%%%%%%%%%%%%%%%%%%%%%%%%%%%%%%%%%%%%%%%%%%%%%
%%%%%%%%%%%%%%%%%%%%%%%%%%%%%%%%%%%%%%%%%%%%%%%%%%%%%%%%%%%%%%%%%%%
%%%%%%%%%%%%%%%%%%%%%%%%%%%%%%%%%%%%%%%%%%%%%%%%%%%%%%%%%%%%%%%%%%%
%%%%%%%%%%%%%%%%%%%    section 2                     %%%%%%%%%%%%%%
\section{Singularity of a fiber outside the branch divisor}

\subsection{Hypertangent divisors and linear systems}

Let $\varphi=\varphi_{F,o}\colon \widetilde F\to F$ be the blow up
of the fiber at an arbitrary point $o$,
$\varphi_G=\varphi_{G,p}\colon \widetilde G\to G$ the blow up of
the fiber $G$ at the point $p=\sigma(o)$, $E=E_F\subset \widetilde
F$ and $E_G\subset \widetilde G$ the exceptional divisors.

{\bf Definition 2.1.} The linear system
$$
\varphi_*(|kH_F-(k+1)E|)
$$
of divisors on $F$ (respectively, the linear system
$$
(\varphi_G)_* (|kH_G-(k+1)E_G|)
$$
of divisors on $G$) is called the $k$-{\it th hypertangent linear
system} and denoted by the symbol $\Lambda_k=\Lambda^F_k$
(respectively, $\Lambda^G_k$).

One can say that $\Lambda_k$ is the largest linear subsystem of
the system $|kH_F|$, the strict transform of which satisfies the
property
$$
{\widetilde \Lambda}_k\subset |kH_F-(k+1)E|,
$$
and similarly for $G$. In the general case one cannot assert that
\begin{equation}
\label{2.1} \sigma^* \Lambda^G_k\subset \Lambda_k,
\end{equation}
since if $p\in W_G$ is a smooth point of the branch divisor, then
the double cover $\sigma\colon F\to G$ {\it does not} extend to a
double cover $\widetilde F\to \widetilde G$ (there is a rational
map of degree two between these varieties; this rational map has a
fairly simple structure, however it is not a finite morphism). But
if $p\not\in W_G$ or $p\in W_G$ is a double point of the branch
divisor, then the inclusion (\ref{2.1}) holds.

The symbol
$$
\Lambda^E_k
$$
stands for the corresponding linear system on the exceptional
divisor:
$$
\Lambda^E_k={\widetilde \Lambda}_k|_E\quad \mbox{or}\quad
\Lambda^E_k={\widetilde \Lambda}_k^G|_{E_G},
$$
depending on the context. It is easy to see that
$$
(\widetilde{\mathop{\rm Bs}\Lambda_k}\circ E)= \mathop{\rm Bs}
\Lambda^E_k
$$
in the scheme-theoretic sense, in particular, the corresponding
effective algebraic cycles are equal, that is, the equality
respects multiplicities.

Abusing out notations, we sometimes use the notion of a
hypertangent system for a certain special subsystem of the
hypertangent system, which permits an explicit description. In
practice it is these special subsystems that we use. Let $p\in G$
be a point, $z_1,\dots,z_{M+1}$ a system of linear coordinates
with the origin at $p$, and assume that the hypersurface $G$ is
given by the equation
$$
f=q_a+q_{a+1}+\dots+q_m=0,
$$
$a=1$ or 2. Then
\begin{equation}
\label{2.2} \Lambda^G_k\supset\left|
\sum^k_{i=a}s_{k-i}f_i\right|,
\end{equation}
where
$$
f_i=q_a+\dots+q_i,
$$
$k\geq a$ and $s_j$ means an arbitrary homogeneous polynomial of
degree $j$ in the variables $z_*$. The inclusion (\ref{2.2}) is
obvious, since
$$
f_i|_G=(-q_{i+1}-\dots-q_m)|_G.
$$
Now assume that $p\not\in W_G$. Let us construct the hypertangent
system $\Lambda_k$. Obviously,
$\Lambda_k\supset\sigma^*\Lambda^G_k$, but in fact the system
$\Lambda_k$ is much larger. Following [19,22,25], let us describe
the construction of hypertangent divisors, associated with the
double cover $\sigma$. Since $p\not\in W_G$, we may assume that
the hypersurface
$$
W_t=W\cap {\mathbb P}_t\subset{\mathbb P}
$$
is given by the equation
$$
g(z_*)=1+w_1+\dots+w_{2l}=0,
$$
$w_i(z_*)$ are homogeneous of degree $i$. Setting formally
\begin{equation}
\label{2.2a}
\sqrt{g}=1+\sum^{\infty}_{i=1}\Phi_i(w_1,\dots,w_{2l}),
\end{equation}
where $\Phi_i(w_1(z_*),\dots,w_{2l}(z_*))$ are homogeneous
polynomials of degree $i$ in $z_*$, write for $j\geq 1$
$$
[\sqrt{g}]_j=1+\sum^j_{i=1}\Phi_i(w_*(z_*)).
$$
Now we get
\begin{equation}
\label{2.3} \Lambda_k\supset\left| \sum^k_{i=a}s_{k-i}f_i+
\sum^{\min\{k,2l-1\}}_{i=l}s^*_{k-i}(y-[\sqrt{g}]_i) \right|,
\end{equation}
where $s^*_{k-i}$ are homogeneous polynomials in $z_*$ of degree
$k-i$; if $k\leq l-1$, then the right-hand side is assumed to be
equal to zero. The inclusion (\ref{2.3}) follows from (\ref{2.2})
and the following fact.

{\bf Lemma 2.1.} {\it In the local coordinates $z_*$ we get
$$
(y-[\sqrt{g}]_i)|_F=2\Phi_{i+1}(w_*(z_*))|_F+\dots,
$$
where the dots stand for a formal series, the components of which
are homogeneous polynomials of degree $i+2$ and higher in the
variables $z_*$.}

{\bf Proof:} it is obvious, since $(y^2-g)|_F\equiv 0$, $g(p)=1$
and the formal decomposition (\ref{2.2a}) holds.

Note that
$$
\Phi_i(w_*)=\frac12 w_i+A_i(w_1,\dots,w_{i-1}).
$$
Now let us consider the case when $p=\sigma(o)\in W_G$. If the
branch divisor is non-singular at the point $p$, then the local
equation of the hypersurface $W_t$ is of the form
$$
g(z_*)=w_1+\dots+w_{2l}=0,
$$
where the linear forms $q_1$, $w_1$ are linearly independent.
Since the inverse image of the divisor
$$
\{w_1|_G=0\}
$$
on $F$ is obviously singular, we obtain:
\begin{equation}
\label{2.5a} \Lambda_k\supset \left| \sum^k_{i=1}s_{k-i}f_i +
s_{k-1}w_1 \right|.
\end{equation}
However, if $p=\sigma(o)$ is a singularity of the divisor $W_G$,
then our methods of constructing hypertangent linear systems give
at most the inclusion $\Lambda_k\supset\sigma^*\Lambda^G_k$, where
$\Lambda^G_k$ is given by the formula (\ref{2.2}).

The regularity conditions make it possible to get a lower bound
for the codimension of the base set of hypertangent systems. In
the formulae below it is assumed that the segment
$[a,b]\subset{\mathbb R}$ is an empty set when $b<a$. For an
arbitrary point $o\in F$ set
$$
{\cal M}=[a,m-1]\cap{\mathbb Z}_+=\{a,\dots,m-1\},
$$
where $a=\mathop{\rm mult}\nolimits_p G\in\{1,2\}$, $p=\sigma(o)$,
and
$$
{\cal L}=[l,2l+a-3]\cap{\mathbb Z}_+=\{l,\dots,2l+a-3\}.
$$
Thus the sets ${\cal M}$, ${\cal L}$ depend on the type of the
point $o$. At each stage of the proof the point $o$ is assumed to
be fixed and the symbols ${\cal M}$, ${\cal L}$ mean the sets
corresponding to this point.

For $e=\max\{m-1,2l-1\}$ we denote the hypertangent linear system
$\Lambda_e$ by the symbol $\Lambda_{\infty}$.

{\bf Proposition 2.1.} {\it The following estimates hold:

{\rm (i)} if $p=\sigma(o)\not\in W_G$ is a smooth point of the
hypersurface $G$, then
$$
\mathop{\rm codim}\nolimits_o \mathop{\rm Bs} \Lambda_k\geq
\mathop{\rm codim}\nolimits_E \mathop{\rm Bs} \Lambda^E_k\geq
\sharp [1,k]\cap{\cal M} + \sharp [1,k]\cap {\cal L},
$$
in particular,
$$
\mathop{\rm dim}\nolimits_o \mathop{\rm Bs} \Lambda_{\infty}\leq
1,
$$

{\rm (ii)} if $p=\sigma(o)\not\in W_G$ is a double point of the
hypersurface $G$, then
$$
\mathop{\rm codim}\nolimits \mathop{\rm Bs} \Lambda_k\geq
\mathop{\rm codim}\nolimits_E \mathop{\rm Bs} \Lambda^E_k\geq
\sharp [2,k]\cap{\cal M} + \sharp [2,k]\cap {\cal L},
$$
and moreover
$$
\mathop{\rm Bs} \Lambda_{\infty}\subset C_*
$$
(see Sec. 1.3). More to that, let $P\subset{\mathbb P}$, $P\ni p$,
be an arbitrary hyperplane, $P_F=\sigma^{-1}(P\cap G)$ the
corresponding section of the fiber $F$,
$\Lambda^P_k=\Lambda_k|_{P_F}$ the restriction of the linear
system $\Lambda_k$ onto $P_F$. Then for $k\leq \max \{m,2l\}-2$
$$
\mathop{\rm codim}\nolimits_{P_F} \mathop{\rm Bs} \Lambda^P_k\geq
\sharp [2,k]\cap{\cal M} + \sharp [2,k]\cap {\cal L},
$$
and
\begin{equation}
\label{2.4} \mathop{\rm dim} \mathop{\rm Bs}
\Lambda^P_{\infty}\leq 1,
\end{equation}
whereas if in (\ref{2.4}) the equality holds then the degree of
the one-dimensional part of the basic subscheme $\mathop{\rm Bs}
\Lambda^P_{\infty}$ does not exceed $\lambda_{m,l}$.

{\rm (iii)} If $p=\sigma(o) \in W_G$ is a smooth point on the
branch divisor $W_G$, then the following inequality holds:
$$
\mathop{\rm codim}\nolimits_o \mathop{\rm Bs} \Lambda_k\geq
\mathop{\rm codim}\nolimits_E \mathop{\rm Bs} \Lambda^E_k\geq
 \sharp [1,k]\cap{\cal M} + 1,
$$

{\rm (iv)} if $p=\sigma(o) \in W_G$ is a double point on the
branch divisor $W_G$, then the following inequality holds:
$$
\mathop{\rm codim}\nolimits_o \mathop{\rm Bs} \Lambda_k\geq
\mathop{\rm codim}\nolimits_E \mathop{\rm Bs} \Lambda^E_k\geq
\sharp [1,k]\cap{\cal M}.
$$
}

{\bf Proof.} To obtain out claims, we replace the hypertangent
linear systems $\Lambda_k$ by their subsystems (\ref{2.2}),
(\ref{2.3}) and (\ref{2.5a}), constructed above, and use the
regularity conditions (Sec. 1.3,1.4). Q.E.D.

%%%%%%%%%%%%%%%%%%%%%%%%%%%%%%%%%%%%%%%%%%%%%%%%%%%%%%%%%%%%%%%%%%%
%%%%%%%%%%%%%%%%%%%%%%%%%%%%%%%%%%%%%%%%%%%%%%%%%%%%%%%%%%%%%%%%%%%
%%%%%%%%%%%%%%%%%%%    subsection 2.2   %%%%%%%%%%%%%%%%%%%%%%%%%%%
\subsection{Scheme of the proof of the condition (vs)}

Assume that there exists a prime divisor $Y\subset F=F_t$,
satisfying the estimate
\begin{equation}
\label{2.5} \frac{\displaystyle \mathop{\rm mult}\nolimits_x
\widetilde Y }{\displaystyle \mathop{\rm deg}\nolimits
Y}>\frac{1}{m},
\end{equation}
where $x\in E$ is an infinitely near point of the first order,
that is, $E\subset \widetilde F$ is the exceptional divisor of the
blow up of the point $o\in F$, $\varphi\colon \widetilde F\to F$.
Here the singular point $o\in F$ is generated by a singularity of
the hypersurface $G=G_t$, that is, $p=\sigma(o)\in G$ is a
non-degenerate double point, $p\not\in W$. Set
$\sigma^{-1}(p)=\{o,o^+\}$ and let $\varphi_G\colon \widetilde
G\to G$ be the blow up of the point $p$. The map $\sigma$ extends
in an obvious way to a morphism
$$
\tilde \sigma\colon \widetilde F\setminus\{o^+\}\to \widetilde G,
$$
whereas on the exceptional divisor $E\subset \widetilde F$ the
morphism $\tilde \sigma$ is an isomorphism, which makes it
possible to identify $E$ with the exceptional divisor of the blow
up $\varphi_G$ and thus consider $E$ as embedded in ${\mathbb T} =
{\mathbb P} (T_p{\mathbb P})\cong {\mathbb P}^M$, that is, in the
exceptional divisor of the blow up
$$
\varphi_{{\mathbb P}}\colon \widetilde {\mathbb P} \to{\mathbb P}
$$
of the point $p\in{\mathbb P}$. Depending on the context one of
the inclusions $E\subset \widetilde F$ or $E\subset \widetilde G$
will be meant.

Let us show that the assumption (\ref{2.5}) leads to a
contradictions. In order to do that, we will use the method
developed in [20]. The arguments break into a few steps. The first
step is given by

{\bf Proposition 2.2.} {\it There exists a hyperplane $P\subset
{\mathbb P}$, $P\ni p$, such that $\sigma(Y)\not\subset P$ and the
effective algebraic cycle $Y_P=(Y\mathop{\circ}\nolimits_F P_F)$,
where $P_F=\sigma^{-1}(P_G)$, $P_G=P\cap G$ is a hyperplane
section, satisfies the estimate
$$
\frac{\mathop{\rm mult}\nolimits_o}{\mathop{\rm deg}} Y_P
>\frac{3}{2m}.
$$
}

The symbol $\mathop{\circ}\nolimits_F$ is used to emphasize that
the cycle $Y_P$ is constructed in the sense of the intersection
theory on $F$, and not on $V$. For a proof of the proposition see
[20].

Step two. Consider the variety $P_F\subset F$. It is an
irreducible variety of dimension $M-1$ with the double point $o\in
P_F$. Let $\varphi_P\colon \widetilde P\to P_F$ be the blow up of
the point $o$, $E_P\subset \widetilde P$ the exceptional divisor.
Obviously, $\widetilde P$ embeds into $\widetilde F$, and $E_P$
into $E$ as a hyperplane section of the quadric $E$ with respect
to the embedding $E\hookrightarrow {\mathbb T}$. Since the variety
$\widetilde F$ is factorial, the strict transform $\widetilde Y$
is a Cartier divisor. Therefore, the effective cycle ${\widetilde
Y}_P=(\widetilde Y \circ \widetilde P)$, that is, the strict
transform of the cycle $Y_P$ on $\widetilde P$, is a Cartier
divisor,
$$
{\widetilde Y}_P \sim aH_P-bE_P,
$$
where $H_P$ is the class of a hyperplane section. By Proposition
2.2
$$
b>\frac32 a.
$$
By the regularity condition we get for the tangent divisor
$$
T=\sigma^{-1}(T_pG\cap G)
$$
that $\mathop{\rm mult}\nolimits_o T=6$, $\mathop{\rm
deg}\nolimits T=4m$, so that for the class of its strict transform
$\widetilde T\subset \widetilde F$ we get $\widetilde T\sim 2H-3E$
and thus for its restriction ${\widetilde T}_P=\widetilde T \cap
\widetilde P$ on $\widetilde P$ we get
$$
{\widetilde T}_P \sim 2H_P-3E_P.
$$

{\bf Proposition 2.3.} {\it Let $Z\sim \alpha H_P-\beta E_P$ be an
effective Cartier divisor on $\widetilde P$. Assume that
$\beta>\frac32 \alpha$. Then $Z$ contains ${\widetilde T}_P $ as a
component of positive multiplicity.}

{\bf Proof} is given below.

Step three. Write down the effective divisor ${\widetilde Y}_P$ in
the following form:
$$
{\widetilde Y}_P=c{\widetilde T}_P+Z,
$$
where $c\in{\mathbb Z}_+$ and the effective divisor $Z$ does not
contain ${\widetilde T}_P$ as a component. Setting $Z\sim \alpha
H_P-\beta E_P$, we obtain from the system of equations
$$
a=2c+\alpha,\quad b=3c+\beta
$$
and the condition $2b>3a$, that
$$
\beta>\frac32 \alpha.
$$
By Proposition 2.3 this implies that ${\widetilde T}_P$ is a
component of positive multiplicity of the divisor $Z$. A
contradiction.

Thus we have proved that the estimate (\ref{2.5}) is impossible
which implies that the condition (vs) holds for the case of a
singular point $o\in F$ outside the branch divisor. Q.E.D.

%%%%%%%%%%%%%%%%%%%%%%%%%%%%%%%%%%%%%%%%%%%%%%%%%%%%%%%%%%%%%%%%%%%
%%%%%%%%%%%%%%%%%%%%%%%%%%%%%%%%%%%%%%%%%%%%%%%%%%%%%%%%%%%%%%%%%%%
%%%%%%%%%%%%%%%%%%%    subsection 2.3   %%%%%%%%%%%%%%%%%%%%%%%%%%%
\subsection{Movable families of curves}

Let us prove Proposition 2.3. We use the method of the paper [20].

{\bf Lemma 2.2.} {\it The divisor $T_P=T\cap P_F$ is swept out by
a family of curves $\{C_{\delta},\delta\in \Delta\}$, the general
member of which is irreducible and satisfies the inequality
\begin{equation}
\label{2.10} \frac{\mathop{\rm mult}\nolimits_o}{\mathop{\rm deg}}
C_{\delta} >\frac23.
\end{equation}
}

First of all, let us obtain Proposition 2.3 from this fact. Let
$\{{\widetilde C}_{\delta},\delta\in \Delta\}$ be the strict
transform of this family of curves on $\widetilde P$, ${\widetilde
T}_P\subset \widetilde P$ the strict transform of the divisor
$T_P$. Obviously,
$$
(Z\cdot {\widetilde C}_{\delta}) =\alpha \mathop{\rm deg}\nolimits
C_{\delta} -\beta \mathop{\rm mult}\nolimits_o C_{\delta} <0,
$$
since $\beta>\frac32 \alpha$. Therefore ${\widetilde
C}_{\delta}\subset Z$. However the curves ${\widetilde
C}_{\delta}$ sweep out ${\widetilde T}_P$, thus $Z\supset
{\widetilde T}_P$. Q.E.D. for Proposition 2.3.

{\bf Proof of Lemma 2.2.} The variety $P_F$ is of dimension
$m+l-2$, the divisor $T_P\subset P_F$ is of dimension $m+l-3$. We
construct the required family of curves
$(C_{\delta},\delta\in\Delta)$, intersecting $T_P$ with $m+l-4$
hypertangent divisors. To order the construction procedure, let us
introduce some new notations:
$$
{\cal M}=\{2,\dots,m-1\},\quad {\cal L}=\{l,\dots,2l-1\},
$$
\begin{equation}
\label{2.6} c_e=\sharp [4,e]\cap{\cal M}+ \sharp [3,e]\cap {\cal
L},\quad e\in {\mathbb Z}_+.
\end{equation}
Here and below we assume silently that the segment $[a,b]\subset
{\mathbb R}$ is the empty set when $b<a$. For $e\leq 2$ we get
$c_e=0$, for $e\geq \max \{m,2l\}-1$ we get that $c_e=m+l-4$,
provided that $l\geq 3$. Let us assume that this is the case and
that $m\geq 4$. Note that for $m=4$
$$
\sharp [4,e]\cap {\cal M}=0,
$$
since this set is empty. The cases $l=2$ and $m=3$ we will treat
separately. Obviously,
$$
c_{e+1}\geq c_e.
$$
Define the {\it ordering function}
$$
\chi\colon \{1,\dots,m+l-4\}\to {\mathbb Z}_+
$$
by the formula
\begin{equation}
\label{2.7} \chi([c_{e-1}+1,c_e]\cap {\mathbb Z}_+)=e.
\end{equation}
In accordance with our remark above, if $c_{e-1}=c_e$, then the
formula (\ref{2.7}) is meaningless, since the set
$[c_{e-1}+1,c_e]$ is empty. Note that
$$
c_{e+1}-c_e\in \{0,1,2\}
$$
by the definition (\ref{2.6}). It is easy to check that
(\ref{2.7}) gives a correct definition of an integer-valued
function $\chi$.

Denote by the symbol $\Lambda^P_i$ the restriction of the
hypertangent system $\Lambda_i$ onto $P_F$. Set
$$
\Lambda^P=\prod^{m+l-4}_{i=1}\Lambda^P_{\chi(i)}.
$$
Note that in this product the hypertangent system $\Lambda_e$ can
appear at most twice, see (\ref{2.6}). Let
$$
{\mathbb D}=\{D_i\in
\Lambda^P_{\chi(i)},\,\,\,i=1,\dots,m+l-4\}\in \Lambda^P
$$
be a general set of hypertangent divisors.

{\bf Definition 2.2.} We say that a family of closed algebraic
sets $(\Gamma_u,u\in U)$ of (co)dimension $i$ on an algebraic
variety $Z$ is a dense movable family if for a general $u\in U$
all irreducible components of the set $\Gamma_u$ are of
(co)dimension $i$ and these components form a family of
irreducible algebraic varieties sweeping out $Z$.

{\bf Lemma 2.3.} {\it For $i=1,\dots,m+l-4$ the closed algebraic
set
$$
R_i({\mathbb D})=\mathop{\bigcap}\limits^i_{j=1} D_i\cap T_P
$$
is for a general ${\mathbb D} \in \Lambda^P$ of codimension $i$ in
$T_P$. For $i=1,\dots,m+l-5$ the family of cycles
$$
(R_i({\mathbb D}),{\mathbb D}\in \Lambda^P)
$$
is a dense movable family of cycles of codimension $i$ on $T$.}

{\bf Proof.} Set $R_0({\mathbb D})=T$ and argue by induction on
$i=1,\dots,m+l-4$. Assume that the claim of the lemma is proved
for $i\leq j\leq m+l-5$ (if $j=0$, then there is nothing to
prove). Set $\chi(j+1)=e$. By definition,
$$
R_{j+1}({\mathbb D})=R_j({\mathbb D})\cap D_{j+1},
$$
where $D_{j+1}\in \Lambda^P_e$ is a general divisor. By definition
of the function $\chi$ we get
$$
j+1 \in [c_{e-1}+1,c_e].
$$
By Proposition 2.1, the following inequality holds:
\begin{equation}
\label{2.8} \mathop{\rm codim}\nolimits_{P_F} \mathop{\rm Bs}
\Lambda^P_e\geq c_e+1,
\end{equation}
so that
$$
\mathop{\rm codim}\nolimits_{T_P} \mathop{\rm Bs}
\Lambda^P_e|_{T_P}\geq c_e,
$$
whereas
$$
\mathop{\rm codim}\nolimits_{T_P} R_j({\mathbb D})=j\leq c_e-1
$$
by (\ref{2.6}), (\ref{2.7}). Therefore, neither of the irreducible
components of the closed subset $R_j({\mathbb D})$ is contained in
the base set of the hypertangent system $\Lambda^P_e$. In
particular,
$$
R_j({\mathbb D})\not\subset D_{j+1}
$$
and therefore $R_{j+1}({\mathbb D})$ is a closed subset of pure
codimension $j+1$ in $T_P$, which proves the first claim of the
lemma. Now assume that $j\leq m+l-6$. Then either
$$
e\leq \max \{m,2l\}-2,
$$
so that by Proposition 2.1 we get the estimate
$$
\mathop{\rm codim}\nolimits_{P_F}\mathop{\rm Bs} \Lambda^P_e\geq
c_e+2,
$$
which is stronger than the inequality (\ref{2.8}), or $e=\max\{
m,2l\}-1$, but in this case $c_e=c_{e-1}+2$, since
$$
j+2\in [c_{e-1}+1,c_e],
$$
so that
$$
\mathop{\rm codim}\nolimits_{T_P}R_j({\mathbb D})=j=c_e-2.
$$
In any case for each irreducible component $Z$ of the set
$R_j({\mathbb D})$ for $j\leq m+l-6$ we get
$$
\mathop{\rm codim}\nolimits_Z \mathop{\rm Bs} (\Lambda^P_e|_Z)
\geq 2,
$$
so that the linear system $\Lambda^P_e|_Z$ is movable. This proves
the second claim of Lemma 2.3.

Consider the family of closed one-dimensional sets
$$
(R({\mathbb D})=R_{m+l-4}({\mathbb D}), {\mathbb D}\in \Lambda^P).
$$
We can no longer assert that irreducible components of the set
$R({\mathbb D})$ form a movable family of curves: at the last
step, that is, when we make curves from surfaces, some fixed
components can appear. However, in any case the following
decomposition holds:
\begin{equation}
\label{2.9} R({\mathbb D})=(T_P\circ D_1 \circ \dots \circ
D_{m+l-4})= \sum_{\delta_i\in\Delta}C_{\delta_i}+\Phi,
\end{equation}
where $(C_{\delta},\delta\in\Delta)$ is a movable family of
curves, $\Phi$ an effective 1-cycle, that is, the fixed part of
the family of curves $R({\mathbb D})$, ${\mathbb D}\in\Lambda^P$.
We get the equality of 1-cycles
$$
\Phi=\mathop{\rm Bs} \Lambda^P_{\infty}.
$$
The family $(C_{\delta},\delta\in\Delta)$ sweeps out $T_F$, {\it
if it is non-empty}. However, by construction
$$
\mathop{\rm deg}\nolimits R({\mathbb D})=
4m\prod^{m+l-4}_{j=1}\chi(j)= 4m\left( \prod^{m-1}_{j=4}j \right)
\left( \prod^{2l-1}_{j=l}j \right)=
$$
$$
=\frac{2m! (2l-1)! }{3 (l-1)!},
$$
whereas by the regularity condition
$$
\mathop{\rm deg}\nolimits \Phi< \lambda_{m,l}=\frac{m! (2l-1)! }{6
(l-1)!}< \mathop{\rm deg}\nolimits R({\mathbb D}).
$$
Therefore, the family of irreducible curves
$(C_{\delta},\delta\in\Delta)$ is non-empty and sweeps out the
divisor $T$.

Let us, finally, estimate the ratio $\mathop{\rm
mult}\nolimits_o/\mathop{\rm deg}\nolimits$ for a general curve
$C_{\delta}$. As we mentioned above, $\mathop{\rm mult}\nolimits_o
\Phi = \mathop{\rm deg}\nolimits \Phi$ (see Sec. 1.3). Besides,
for a general set ${\mathbb D}\in\Lambda^P$ the ratio
$$
\frac{\mathop{\rm mult}\nolimits_o}{\mathop{\rm deg}} C_{\delta_i}
$$
(in the sense of the formula (\ref{2.9})) does not depend on $i$.
Consequently,
$$
\frac{\mathop{\rm mult}\nolimits_o}{\mathop{\rm deg}} C_{\delta}=
\frac{\mathop{\rm mult}\nolimits_o R({\mathbb D})-\mathop{\rm
deg}\nolimits \Phi }{\mathop{\rm deg}\nolimits R({\mathbb
D})-\mathop{\rm deg}\nolimits \Phi }.
$$
However, by construction
$$
\mathop{\rm mult}\nolimits_o R({\mathbb D})\geq
6\prod^{m+l-4}_{j=1}(\chi(j)+1)= 6 \left( \prod^{m}_{j=5}j \right)
\left( \prod^{2l}_{j=l+1}j \right)=
$$
$$
=\frac{m!}{4}\cdot \frac{(2l)!}{l!},
$$
whence
$$
\frac{\mathop{\rm mult}\nolimits_o}{\mathop{\rm deg}}
C_{\delta}\geq \frac{\displaystyle\frac{m!}{4}\cdot
\frac{(2l)!}{l!}-\lambda_{m,l} }{\displaystyle\frac{2m!}{3}\cdot
\frac{(2l-1)!}{(l-1)!}-\lambda_{m,l} }>\frac23
$$
in accordance with the choice of the number $\lambda_{m,l}$. This
proves the lemma for $m\geq 5$, $l\geq 3$.

If $l=2$, then the arguments presented above work with the only
modification: instead of (\ref{2.6}) one should use the formula
$$
c_e=\sharp [3,e]\cap {\cal M} +\sharp [3,e]\cap {\cal L},
$$
$e\in{\mathbb Z}_+$. In this case an independent hypertangent
divisor adds into the linear system $\Lambda^P_2$ and thus the
codimension of its base set (and the codimension of the base set
of all the subsequent hypertangent systems $\Lambda^P_j$, $j\geq
3$) exceeds by one the corresponding codimension in the just
considered case $l\geq 3$. It is this fact that makes it possible
to change the definition of the number $c_e$ and accordingly shift
by one the function $\chi$. The rest of the arguments are
completely similar to the case $l\geq 3$ discussed above.

The case $m\geq 3$ is slightly harder. In order to obtain the
needed codimension of the base set of a hypertangent system, one
should use the following set of hypertangent divisors:
$$
{\mathbb D}=\{ D_i\in\Lambda^P_{l+i}\,\, |\,\, i=1,\dots,l-1\}.
$$
We draw the reader's attention to the fact that the first divisor
in this set is taken from the linear system $\Lambda^P_{l+1}$,
that is, in contrast to the case $m\geq 4$, which we considered
above, we skip the system $\Lambda^P_l$. As a result, we obtain
once again a movable family of closed algebraic sets
$$
R_k({\mathbb D})=
\left(\mathop{\bigcap}\limits^k_{j=1}D_j\right)\cap T
$$
for $k\leq 2l-2$, whereas irreducible components of the sets
$R_k({\mathbb D})$ form a family and sweep out $T$. Again we
modify the family of curves $R_{2l-1}({\mathbb D})$, deleting the
fixed part $\Phi$ of degree $\mathop{\rm deg}\nolimits
\Phi<\lambda_{3,l}= 12(l-2)\frac{\displaystyle
(2l-1)!}{\displaystyle (l+1)!}$ and obtain a family of irreducible
curves $(C_{\delta},\delta\in\Delta)$, sweeping out $T$ and
satisfying the estimate (\ref{2.10}). Proof of Lemma 2.2 is now
complete.

%%%%%%%%%%%%%%%%%%%%%%%%%%%%%%%%%%%%%%%%%%%%%%%%%%%%%%%%%%%%%%%%%%%
%%%%%%%%%%%%%%%%%%%%%%%%%%%%%%%%%%%%%%%%%%%%%%%%%%%%%%%%%%%%%%%%%%%
%%%%%%%%%%%%%%%%%%%%%%%%%%%%%%%%%%%%%%%%%%%%%%%%%%%%%%%%%%%%%%%%%%%
%%%%%%%%%%%%%%%%%%%%%%%%%%%%%%%%%%%%%%%%%%%%%%%%%%%%%%%%%%%%%%%%%%%
%%%%%%%%%%%%%%%%%%%    section 3                      %%%%%%%%%%%%%
\section{Singularity of a fiber on the branch divisor}
\subsection{Notations and discussion of the regularity condition}

We have the double cover
$$
F=F_t\stackrel{\sigma}{\to}G= G_t\subset {\mathbb P}={\mathbb
P}^{M+1},
$$
$G\subset{\mathbb P}$ is a smooth hypersurface of degree $m\leq
M-1$. At the point $p\in G$ the branch hypersurface $W_G=W\cap G$
has an isolated quadratic singularity, so that
$o=\sigma^{-1}(p)\in F$ is an (isolated) non-degenerate double
point of the fiber $F$. We get the following commutative diagram
of maps
$$
\begin{array}{rrcccll}
E=E_F & \subset & \widetilde F & \stackrel{\tilde \sigma}{\to}
& \widetilde G & \supset & E_G \\
    & \varphi_F & \downarrow &   &  \downarrow & \varphi_G &   \\
    &        &  F & \stackrel{\sigma}{\to} & G, &  &
\end{array}
$$
where $\varphi_F$ and $\varphi_G$ are the blow ups of the points
$o\in F$ and $p\in G$, respectively, $E_F$ and $E_G$ are the
exceptional divisors, $\tilde \sigma$ the double cover, branched
over ${\widetilde W}_G \subset \widetilde G$, that is, over the
strict transform of the divisor $W_G$. Besides,
$$
{\tilde \sigma}_E={\tilde \sigma}|_E\colon E \to E_G\cong {\mathbb
P}^{M-1}
$$
is the double cover, branched over the quadric
$$
W_E={\widetilde W}_G\cap E_G.
$$
The symbol $H_E$ stands for the hyperplane section of the quadric
$E$ with respect to the standard embedding $E\hookrightarrow
{\mathbb P}^M$, $\mathop{\rm Pic} E={\mathbb Z} H_E$.

Let $W_t=W\cap {\mathbb P}_t$ be given by the equation
$$
h=w_1+w_2+\dots +w_{2l}=0,
$$
and $G$ by the equation
$$
f=q_1+q_2+\dots +q_m=0
$$
with respect to the affine coordinates $z_*=(z_1,\dots,z_{M+1})$
with the origin at the point $p$. The divisor $W_G$ has at $p$ a
non-degenerate quadratic singularity, so that $w_1=\lambda q_1$,
for simplicity of notations assume that $q_1=z_{M+1}$. The
quadratic polynomial ${\bar w}_2=w_2|_{\{z_{M+1}=0\}}$ is of the
maximal rank. Take $z_1,\dots,z_{M}$ for homogeneous coordinates
on $E_G$, then
$$
{\tilde \sigma}_E\colon E\to E_G\cong {\mathbb P}^{M-1}
$$
is branched over the non-singular quadric $W_E=\{{\bar w}_2=0\}$.
For an arbitrary point $y\in E_G\setminus W_E$ let $C(y)\subset
E_G$ be the cone consisting of all lines $L\subset E_G$ that
contain $y$ and touch $W_E$. More formally, let
$$
\pi_y\colon E_G\setminus \{y\} \to {\mathbb P}^{M-2}
$$
be the projection from the point $y$. Its restriction onto the
quadric $W_E$,
$$
\pi_y|_{W_E}\colon W_E \to {\mathbb P}^{M-2}
$$
is a double cover, branched over a quadric $Q(y)\subset {\mathbb
P}^{M-2}$. Now
$$
C(y)=\overline{\pi^{-1}_y(Q(y)) }.
$$
Obviously, $C(y)$ is a quadric cone with the vertex at the point
$y$. Since the quadric $W_E$ is non-singular, the cone $C(y)$ has
only one singularity, that is, the point $y$.

Denote the restriction of the polynomial $q_i$ onto the hyperplane
$q_{M+1}=0$ by the symbol ${\bar q}_i$. By the regularity
condition, the system of homogeneous equations
$$
{\bar q}_2=\dots={\bar q}_m=0
$$
defines in $E_G$ an irreducible reduced complete intersection of
codimension $(m-1)$, an irreducible subvariety $Z_{2\cdot \dots
\cdot m}$. Moreover, the quadric ${\bar q}_2=0$ is smooth and
distinct from $W_E$.

{\bf Lemma 3.1.} {\it Assume that the condition {\rm (R2.2)}
holds. Then the subvariety $Z_{2\cdot \dots \cdot m}$ is not
contained in a quadric cone $C(y)$, $y\in E_G\setminus W_E$, and
in a tangent plane $T_y W_E$, $y\in W_E$.}

{\bf Proof.} Set
$$
Z_{2\cdot \dots \cdot j}= \{ z\in {\mathbb P}^{M-1}\,\, |\,\,
{\bar q}_2=\dots={\bar q}_j=0\}.
$$
It is easy to see that $Z_{2\cdot \dots \cdot j}$ is an
irreducible reduced complete intersection of codimension $j$. From
the long exact cohomology sequence we obtain that
$$
h^0({\cal O}_{Z_{2}}(2))=\dots = h^0({\cal O}_{Z_{2\cdot \dots
\cdot j}}(2))=\dots = h^0({\cal O}_{Z_{2\cdot \dots \cdot m}}(2)),
$$
and moreover, the restriction map
$$
H^0({\cal O}_{{\mathbb P}^{M-1}}(2))\to H^0({\cal O}_{Z_{2\cdot
\dots \cdot m}}(2))
$$
is surjective. This implies that $Z_{2\cdot \dots \cdot m}$ is
contained in one and only one quadric $Z_2$ and thus is not
contained in any quadric cone $C(y)$, $y\in E_G\setminus W_E$. In
a similar way, the restriction map
$$
H^0({\cal O}_{{\mathbb P}^{M-1}}(1))\to H^0({\cal O}_{Z_{2\cdot
\dots \cdot m}}(1))
$$
is an isomorphism, so that $Z_{2\cdot \dots \cdot m}$ is not
contained in a hyperplane, in particular, in a hyperplane of the
form $T_y W_E$, $y\in W_E$.

Now fix a prime divisor $R\subset F$, and let $\widetilde R
\subset \widetilde F$ be its strict transform. Fix also an
arbitrary point $x\in E$, lying outside the branch divisor of the
cover ${\tilde \sigma}_E$, that is, ${\tilde \sigma}(x)\not\in
W_E$. (For a point $x\in E$ on the branch divisor the arguments
given below work automatically with simplifications. The arguments
of Sec. 2 can be also used in this case, in contrast to the
situation outside the branch divisor $W_E$.)

{\bf Proposition 3.1.} {\it The following estimate holds:}
\begin{equation}
\label{3.1} \mu=\mathop{\rm mult}\nolimits_x \widetilde R \leq
\frac{1}{m} \mathop{\rm deg}\nolimits R.
\end{equation}

{\bf Remark.} For some $k\geq 1$ we have $R\sim H_F$, where
$H_F=\sigma^* H_G$ is a hyperplane section. Since obviously
$\mathop{\rm deg}\nolimits R=2mk$, the estimate (\ref{3.1}) takes
the form of the following inequality:
$$
\mu\leq 2k.
$$

%%%%%%%%%%%%%%%%%%%%%%%%%%%%%%%%%%%%%%%%%%%%%%%%%%%%%%%%%%%%%%%%%%%
%%%%%%%%%%%%%%%%%%%%%%%%%%%%%%%%%%%%%%%%%%%%%%%%%%%%%%%%%%%%%%%%%%%
%%%%%%%%%%%%%%%%%%    subsection 3.2   %%%%%%%%%%%%%%%%%%%%%%%%%%%%
\subsection{Start of the proof of the condition (vs)}

Assume the converse: $\mu> 2k$. We have the presentation
$$
\widetilde R\sim k\varphi^*_F H_F-\nu E,
$$
whereas $\mathop{\rm mult}\nolimits_o R=2\nu$.

{\bf Lemma 3.2.} {\it The following inequality holds: $\nu\leq
2k$.}

{\bf Proof.} Assume the converse: $\nu>2k$. Then
$$
\frac{\mathop{\rm mult}\nolimits_o}{\mathop{\rm deg}} R
>\frac{2}{m}.
$$
Set ${\bar R}=\sigma(R)\subset G$. It is a prime divisor on the
smooth hypersurface $G\subset{\mathbb P}$. Since $\sigma\colon
R\to{\bar R}$ is a finite morphism, we get the inequality
$$
\frac{\mathop{\rm mult}\nolimits_o}{\mathop{\rm deg}} {\bar R}
>\frac{2}{m}.
$$
However, this is impossible, since $p\in G$ is a regular point.
Indeed, the tangent divisor $T_1^+=T_pG\cap G$ satisfies the
equality
$$
\frac{\mathop{\rm mult}\nolimits_p}{\mathop{\rm deg}}
T^+_1=\frac{2}{m},
$$
so that ${\bar R}\neq T^+_1$ and $({\bar R}\circ T^+_1)$ is an
effective cycle of codimension two ($T^+_1$ is obviously
irreducible).

Since $\mathop{\rm mult}\nolimits_p T^+_1=2$, we get the
inequality
$$
\mathop{\rm mult}\nolimits_p ({\bar R}\circ T^+_1) \geq
2\mathop{\rm mult}\nolimits_p {\bar R}.
$$
Taking into account that $\mathop{\rm deg}\nolimits ({\bar R}\circ
T^+_1) =\mathop{\rm deg}\nolimits {\bar R}$, we conclude that
there exists an irreducible component $Y_2$ of the cycle $({\bar
R}\circ T^+_1)$, satisfying the estimate
$$
\frac{\mathop{\rm mult}\nolimits_p}{\mathop{\rm deg}} Y_2 \geq
2\frac{\mathop{\rm mult}\nolimits_p}{\mathop{\rm deg}} {\bar R}.
$$
As usual, let $f=q_1+q_2+\dots+q_m$ be the polynomial defining the
hypersurface $G$ with respect to the coordinate system $z_*$ with
the origin at the point $p$. Setting
$$
f=q_1+q_2+\dots+q_i,
$$
let us construct the hypertangent systems
$$
\Lambda^G_i=\left| \sum^i_{j=1} f_j s_{i-j}|_G=0 \right|,
$$
and consider the standard hypertangent divisors
$$
T^+_i=\{f_i|_G=0\}\in \Lambda^G_i.
$$
Set
$$
T_i=\sigma^* T^+_i, \quad \Lambda_i=\sigma^* \Lambda^G_i.
$$
These divisors and linear systems will be of crucial importance
below. At the moment, note that by the regularity condition we get
$$
\mathop{\rm codim}\nolimits_G \mathop{\rm Bs} \Lambda^G_i=i
$$
(in fact, $\mathop{\rm Bs} \Lambda^G_i=T^+_1\cap \dots \cap
T^+_i$). Let
$$
{\mathbb D}=(D_1,\dots,D_{m-1})\in \prod^{m-1}_{j=1}\Lambda^G_i
$$
be a general set of divisors. Let us construct by induction a
sequence of irreducible subvarieties $Y_i$, $i=1,\dots,m-1$,
satisfying the following properties:

(i) $Y_1={\bar R}$, $Y_2$ was constructed above, $\mathop{\rm
codim}\nolimits_G Y_i=i$;

(ii) $Y_{i+1}\subset Y_i$, $Y_i\not\subset D_{i+1}$, $Y_{i+1}$ is
an irreducible component of the closed set $Y_i\cap D_{i+1}$;

(iii) the estimate
$$
\frac{\mathop{\rm mult}\nolimits_p}{\mathop{\rm deg}} Y_{i+1}\geq
\frac{i+2 }{i+1}\cdot \frac{\mathop{\rm
mult}\nolimits_p}{\mathop{\rm deg}} Y_i
$$
holds.

It is possible to construct this sequence because
$$
\mathop{\rm codim}\nolimits_G \Lambda^G_{i+1}=i+1
> \mathop{\rm codim}\nolimits_G Y_i,
$$
so that for a general divisor $D_{i+1}\subset \Lambda^G_{i+1}$ we
have $Y_i\not\subset D_{i+1}$. One can ensure that the property
(iii) holds since $\Lambda^G_j\subset |jH_G|$ and $\mathop{\rm
mult}\nolimits_p \Lambda^G_j=j+1$.

Now for an irreducible subvariety $Y=Y_{m-1}$ we get the estimate
$$
1\geq \frac{\mathop{\rm mult}\nolimits_p}{\mathop{\rm deg}} Y\geq
\underbrace{\frac{m}{m-1}\cdot\frac{m-1}{m-2}\cdot\cdots\cdot
\frac43 \cdot \frac21 }_{
\begin{array}{c}
\| \\ \displaystyle\frac{2m}{3}
\end{array}
}\cdot \frac{\mathop{\rm mult}\nolimits_p}{\mathop{\rm deg}} {\bar
R},
$$
whence we get
$$
\frac{\mathop{\rm mult}\nolimits_p}{\mathop{\rm deg}} {\bar R}
\leq \frac{3}{2m}.
$$
Therefore the ratio $\mathop{\rm mult}\nolimits_p/\mathop{\rm
deg}\nolimits$ attains its maximum at the tangent divisor $T_p
G\cap G$ and this maximum is equal to 2. A contradiction. Q.E.D.
for the lemma.

%%%%%%%%%%%%%%%%%%%%%%%%%%%%%%%%%%%%%%%%%%%%%%%%%%%%%%%%%%%%%%%%%%%
%%%%%%%%%%%%%%%%%%%%%%%%%%%%%%%%%%%%%%%%%%%%%%%%%%%%%%%%%%%%%%%%%%%
%%%%%%%%%%%%%%%%%%    subsection 3.3   %%%%%%%%%%%%%%%%%%%%%%%%%%%%
\subsection{Hypertangent divisors and tangent cones}

Thus $\nu\leq 2k<\mu$. On the other side,
$$
\mu=\mathop{\rm mult}\nolimits_x \widetilde R \leq \mathop{\rm
mult}\nolimits_x(\widetilde R \circ E)\leq \mathop{\rm
deg}\nolimits (\widetilde R \circ E)=2\nu.
$$
Set $B=T_xE\cap E$, where the quadric $E$ is considered as
embedded in ${\mathbb P}^M$ in the standard way. By Lemma 5 from
Sec. 3.5 in [20],
$$
\mathop{\rm mult}\nolimits_B \widetilde R\geq \frac12 (\mu-\nu),
$$
whereas for the effective cycle $R_E=(\widetilde R \circ E)$ we
get
\begin{equation}
\label{3.77} \mathop{\rm mult}\nolimits_B R_E\geq \mu-\nu.
\end{equation}
Set ${\widetilde T}_i, {\widetilde \Lambda}_i$ to be the strict
transforms of the divisors $T_i$ and linear systems $\Lambda_i$ on
the blow up $\widetilde F$ of the fiber $F$. It is easy to see
that ${\widetilde T}_i \subset {\widetilde \Lambda}_i$, since by
construction
$$
\mathop{\rm mult}\nolimits_o \Lambda_i=\mathop{\rm
mult}\nolimits_o T_i.
$$
Set also
$$
{\mathbb T}_i=({\widetilde T}_i \circ E)={\widetilde T}_i\cap E
$$
to be the projectivized tangent cone to the divisor $T_i$ at the
point $o$. Recall that the quadric $E$ is realized as the double
cover ${\tilde \sigma}_E\colon E\to E_G\cong {\mathbb P}^{M-1}$,
branched over the quadric $W_E$. For a system
$(z_1,\dots,z_{M+1})$ of affine coordinates on ${\mathbb P}$ with
the origin at the point $p$ we may assume that $q_1\equiv z_{M+1}$
and therefore $(z_1,\dots,z_{M})$ can be taken for homogeneous
coordinates on the projective space $E_G$. In terms of these
coordinates the hypersurface ${\mathbb T}_i\subset E$ is given by
the equation
$$
({\tilde \sigma}_E)^* q_{i+1}|_{E_G}.
$$
Finally, set
$$
\Lambda^E_i={\widetilde \Lambda}_i|_E
$$
to be the projectivized tangent system of the linear system
$\Lambda_i$ at the point $o$. Equations of divisors of this linear
system are obtained by pulling back to $E$ via ${\tilde \sigma}_E$
the equations
\begin{equation}
\label{3.2} \sum^i_{j=1}{\bar q}_{j+1}{\bar s}_{i-j},
\end{equation}
where $\bar\sharp$ means the restriction of the polynomial
$\sharp$ onto the hyperplane $z_{M+1}=0$. Obviously,
$$
{\mathbb T}_i\sim (i+1) H_E, \quad \Lambda^E_i\subset |(i+1)H_E|,
$$
besides the equations (\ref{3.2}) imply directly that
$$
\mathop{\rm Bs} \Lambda_i=T_1\cap\dots\cap T_i, \quad \mathop{\rm
Bs} \Lambda^E_i= {\mathbb T}_1\cap\dots\cap {\mathbb T}_i,
$$
both equalities in the scheme-theoretic sense.

%%%%%%%%%%%%%%%%%%%%%%%%%%%%%%%%%%%%%%%%%%%%%%%%%%%%%%%%%%%%%%%%%%%
%%%%%%%%%%%%%%%%%%%%%%%%%%%%%%%%%%%%%%%%%%%%%%%%%%%%%%%%%%%%%%%%%%%
%%%%%%%%%%%%%%%%%%    subsection 3.4   %%%%%%%%%%%%%%%%%%%%%%%%%%%%
\subsection{Constructing new cycles}

By the regularity condition the set ${\mathbb T}_1\cap \dots\cap
{\mathbb T}_i$ is irreducible and not contained in the divisor $B$
for all $i=1,\dots,m-1$. Let
$$
{\cal L}= (L_2,\dots,L_{m-1})\in \Lambda_2\times\dots\times
\Lambda_{m-1}
$$
be a general set of hypertangent divisors. We denote the strict
transform of the hypertangent divisor $L_j$ on $\widetilde F$ by
the symbol ${\widetilde L}_j$ and its projectivized  tangent cone
at the point $o\in F$ by the symbol
$$
{\mathbb L}_j=({\widetilde L}_j\circ E).
$$
For a general divisor $L_j\in\Lambda_j$ we get ${\mathbb
L}_j={\widetilde L}_j\cap E$.

{\bf Lemma 3.3.} {\it {\rm (i)} Let $Y\subset F$ be a fixed
irreducible subvariety of codimension $l\leq m-2$. For a general
divisor $L_{l+1}\in \Lambda_{l+1}$ we have $Y\not\subset L_{l+1}$.

{\rm (ii)} Let $Y\subset E$ be a fixed irreducible subvariety of
codimension $l\leq m-2$. For a general divisor $L_{l+1}\in
\Lambda_{l+1}$ we have $Y\not\subset {\mathbb L}_{l+1}$.}

{\bf Proof.} By the regularity condition
$$
\mathop{\rm codim}\nolimits_F \mathop{\rm Bs}\Lambda_{l+1}=l+1,
\quad \mathop{\rm codim}\nolimits_E \mathop{\rm
Bs}\Lambda^E_{l+1}=l+1
$$
and for a general divisor $L_j\in\Lambda_j$ we have ${\mathbb
L}_j\in\Lambda^E_j$. Q.E.D. for the lemma.

{\bf Corollary 3.1.} {\it For a general set ${\cal L}$ we have}
$$
\mathop{\rm codim}\nolimits_F (R\cap L_2\cap\dots\cap
L_{m-1})=m-1,
$$
$$
\mathop{\rm codim}\nolimits_E (R_E\cap {\mathbb L}_2\cap\dots\cap
{\mathbb L}_{m-1})=m-1.
$$

From this fact we obtain that the following effective algebraic
cycles of codimension $m-1$ are well defined on $F$ and $E$,
respectively:
$$
R^+=(R\circ L_2\circ\dots\circ L_{m-1})
$$
and
$$
R^+_E=(R_E\circ {\mathbb L}_2\circ\dots\circ {\mathbb L}_{m-1}),
$$
whereas (for a general set ${\cal L}$)
$$
R^+_E=({\widetilde R}^+\circ E)
$$
is the projectivized tangent cone to the cycle $R^+$ at the point
$o$. Let us describe the structure of these effective cycles.
First of all we get
$$
\mathop{\rm deg}\nolimits R^+=2km\cdot (m-1)!=2km!,
$$
$$
\mathop{\rm mult}\nolimits_o R^+ =\mathop{\rm deg}\nolimits R^+_E=
2\nu\cdot 3\cdot\dots\cdot m=\nu m!.
$$

{\bf Lemma 3.4.} {\it Let $Y$ be an irreducible component of the
cycle $R^+$. If $Y\subset T_1$, then}
$$
Y=T_1\cap T_2\cap\dots\cap T_{m-1}.
$$

{\bf Proof.} By construction, the equation of the divisor $L_i$ is
of the form
$$
f_1s_{i-1}+f_2s_{i-2}+\dots+f_is_0,
$$
where $s_j$ is a homogeneous polynomial of degree $j$ in the
coordinates $z_*$. Since the hypertangent divisors $L_i$ are
assumed to be general, we may assume that $s_0\neq 0$ and thus
normalize the equation by the condition that $s_0=1$. Assume that
$Y\subset T_1$. Then the following polynomials vanish on $Y$:
$$
\begin{array}{lclcccccc}
f_1, & & & & & & & &  \\
f_1s_{2,1} & + & f_2, & & & & & & \\
f_1s_{3,2} & + & f_2 s_{3,1} & + & f_3, & & & &  \\
& & & & \dots & & & & \\
f_1 s_{m-1,m-2} & + &  & \dots & & + & f_{m-2}s_{m-1,1} & + &
f_{m-1},
\end{array}
$$
where $s_{i,j}$ is a homogeneous polynomial of degree $j$. Thus
$$
f_1|_Y\equiv f_2|_Y\equiv \dots \equiv f_{m-1}|_Y\equiv 0,
$$
so that $Y\subset T_1\cap T_2\cap\dots\cap T_{m-1}$, but the
latter set is irreducible and of the same dimension as $Y$. This
proves Lemma 3.4.

%%%%%%%%%%%%%%%%%%%%%%%%%%%%%%%%%%%%%%%%%%%%%%%%%%%%%%%%%%%%%%%%%%%
%%%%%%%%%%%%%%%%%%%%%%%%%%%%%%%%%%%%%%%%%%%%%%%%%%%%%%%%%%%%%%%%%%%
%%%%%%%%%%%%%%%%%%    subsection 3.5   %%%%%%%%%%%%%%%%%%%%%%%%%%%%
\subsection{Degrees and multiplicities}

Set
$$
T=T_1\cap T_2\cap\dots\cap T_{m-1}, \quad {\mathbb T}={\mathbb
T}_1\cap {\mathbb T}_2\cap \dots\cap {\mathbb T}_{m-1}.
$$
Taking into consideration that $T=(T_1\circ \dots \circ T_{m-1})$
and ${\mathbb T}=({\mathbb T}_1\circ \dots \circ {\mathbb
T}_{m-1})_E$, it is easy to verify that
$$
\mathop{\rm deg}\nolimits T=\mathop{\rm mult}\nolimits_o T
=\mathop{\rm deg}\nolimits {\mathbb T}=2m!.
$$
Now write down
\begin{equation}
\label{3.3} R^+=aT+R^{\sharp},\quad R^+_E=a{\mathbb T}
+R^{\sharp}_E,
\end{equation}
where $a\in{\mathbb Z}_+$, the effective cycle $R^{\sharp}$ is
uniquely defined by the condition that it does not contain the
subvariety $T$ as a component, and
$$
R^{\sharp}_E=({\widetilde R}^{\sharp}\circ E)
$$
is the projectivized tangent cone to $R^{\sharp}$ at the point
$o$. Note that the irreducible subvariety ${\mathbb T}$, generally
speaking, can come into the effective cycle $R^{\sharp}_E$ as a
component.

{\bf Lemma 3.5.} {\it The following estimate holds:}
$$
2\mathop{\rm mult}\nolimits_o R^{\sharp}\leq \mathop{\rm
deg}\nolimits R^{\sharp}.
$$

{\bf Proof.} Let $Y$ be an irreducible component of the cycle
$R^{\sharp}$. By construction, $Y\neq T$; therefore by Lemma 3.4
$Y\not\subset T_1$. Thus the closed subset
$$
T_1\cap \mathop{\rm Supp}R^{\sharp}
$$
is of codimension $m$, so that the effective cycle
$$
R^*=(R^{\sharp}\circ T_1)
$$
is well defined. Now we have a standard chain of estimates:
$$
2\mathop{\rm mult}\nolimits_o R^{\sharp}\leq \mathop{\rm
mult}\nolimits_o R^*\leq \mathop{\rm deg}\nolimits R^*=
\mathop{\rm deg}\nolimits R^{\sharp},
$$
which is what we need.

As in Corollary 3.1, Lemma 3.5 implies that the set
$$
B\cap{\mathbb L}_2\cap \dots \cap{\mathbb L}_{m-1}
$$
is of codimension $m-1$ in $E$. Denote by $B^+$ the part of the
effective equidimensional cycle $R^+_E$, the support of which is
contained in $B$:
$$
R^+_E=\sum_{i\in I}r_i Y_i,\quad B^+= \sum_{i\in I, Y_i\subset B}
r_i Y_i.
$$

{\bf Lemma 3.6.} {\it The following estimate holds:}
$$
\mathop{\rm deg}\nolimits B^+\geq (\mu-\nu)m!
$$

{\bf Proof.} Indeed, by (\ref{3.77}) we get
$$
R_E=(\mu-\nu)B+\Delta,
$$
where $\Delta$ is an effective cycle. Furthermore,
$$
\mathop{\rm deg}\nolimits (B\circ {\mathbb L}_2\circ
\dots\circ{\mathbb L}_{m-1})=2\cdot 3\cdot\dots\cdot m=m!,
$$
which proves the lemma.

{\bf Lemma 3.7.} {\it Let $Y\subset {\mathbb
L}_2\cap\dots\cap{\mathbb L}_{m-1}$ be an irreducible subvariety
of codimension $m-1$ in $E$. If $Y\subset {\mathbb T}_1$, then
$Y={\mathbb T}$. }

{\bf Proof.} The equation of the divisor ${\mathbb L}_i$ with
respect to the homogeneous coordinates $z_*$ is of the form
$$
q_2s_{i-1}+\dots+q_{i+1},
$$
where $s_j$ is a homogeneous polynomial of degree $j$. If
$Y\subset {\mathbb T}_1$, then the following polynomials vanish on
$Y$:
$$
\begin{array}{lcccccccc}
q_2, & & & & & & & &  \\
q_2s_{2,1} & + & q_3, & & & & & & \\
& & & & \dots & & & & \\
q_2 s_{m-1,m-2} & + &  & \dots & & + & q_{m-1}s_{m-1,1} & + &
q_{m},
\end{array}
$$
where $\mathop{\rm deg}\nolimits s_{i,j}=j$. Consequently,
$$
q_2|_Y\equiv q_3|_Y\equiv \dots \equiv q_{m}|_Y\equiv 0,
$$
that is, $Y\subset {\mathbb T}$, and since the dimensions
coincide, $Y={\mathbb T}$. Q.E.D. for the lemma.

{\bf Corollary 3.2.} {\it None of the components of the closed set
$$
B\cap {\mathbb L}_2\cap \dots\cap {\mathbb L}_{m-1}
$$
is contained in ${\mathbb T}_1$.}

{\bf Proof.} Let $Y$ be such component and $Y\subset {\mathbb
T}_1$. By the previous lemma, $Y={\mathbb T}$. Thus ${\mathbb
T}\subset B$: a contradiction with the regularity condition.
Q.E.D. for the corollary.

Let us complete, at long last, the proof of Proposition 3.1. From
the presentations (\ref{3.3}) we get
$$
\mathop{\rm deg}\nolimits R^+=2km!= 2am!+\mathop{\rm deg}\nolimits
R^{\sharp},
$$
$$
\mathop{\rm mult}\nolimits_o R^+=\nu m!= 2am!+\mathop{\rm
mult}\nolimits_o  R^{\sharp}.
$$
By Corollary 3.2 the effective cycle $B^+$ lies entirely in
$R^{\sharp}_E$. In particular,
\begin{equation}
\label{3.4} \mathop{\rm deg}\nolimits R^{\sharp}_E\geq \mathop{\rm
deg}\nolimits B^+\geq (\mu-\nu)m!.
\end{equation}
However, $\mathop{\rm deg}\nolimits R^{\sharp}_E=\mathop{\rm
mult}\nolimits_o R^{\sharp}$. Applying Lemma 3.5, we obtain:
$$
2(\nu m!-2a m!)\leq 2k m!-2a m!.
$$
Let us rewrite the inequality (\ref{3.4}) in the form
$$
\nu m!-2a m!\geq (\mu-\nu)m!.
$$
Easy computations give us the two inequalities
$$
k+a\geq \nu,
$$
$$
2\nu-2a\geq \mu,
$$
which imply the desired estimate (\ref{3.1}) in an obvious way.

However, we assumed that $\mu> 2k$. The contradiction completes
our proof of Proposition 3.1 and Theorem 1 as well.

%%%%%%%%%%%%%%%%%%%%%%%%%%%%%%%%%%%%%%%%%%%%%%%%%%%%%%%%%%%%%%%%%%%
%%%%%%%%%%%%%%%%%%%%%%%%%%%%%%%%%%%%%%%%%%%%%%%%%%%%%%%%%%%%%%%%%%%
%%%%%%%%%%%%%%%%%%%%%%%%%%%%%%%%%%%%%%%%%%%%%%%%%%%%%%%%%%%%%%%%%%%
%%%%%%%%%%%%%%%%%%%%%%%%%%%%%%%%%%%%%%%%%%%%%%%%%%%%%%%%%%%%%%%%%%%
%%%%%%%%%%%%%%%%%%    references                       %%%%%%%%%%%%
\section*{References}

{\small

\noindent 1. Algebraic surfaces. By the members of the seminar of
I.R.Shafarevich. I.R.Shafarevich ed. Proc. Steklov Math. Inst.
{\bf 75}. 1965. English transl. by AMS, 1965. 281 p.
\vspace{0.3cm}

\noindent 2. Brown G., Corti A. and Zucconi F. Birational geometry
of 3-fold Mori fibre spaces. Preprint, 2003, 40 p. arXiv:
math.AG/0307301. \vspace{0.3cm}

\noindent 3. Corti A., Pukhlikov A. and Reid M., Fano 3-fold
hypersurfaces, in ``Explicit Birational Geometry of Threefolds'',
London Mathematical Society Lecture Note Series {\bf 281} (2000),
Cambridge University Press, 175-258. \vspace{0.3cm}

\noindent 4. Corti A., Factoring birational maps of threefolds
after Sarkisov. J. Algebraic Geom. {\bf 4} (1995), no. 2, 223-254.
\vspace{0.3cm}

\noindent 5. Graber T., Harris J. and Starr J. Families of
rationally connected varieties. J. Amer. Math. Soc. {\bf 16}
(2002), no. 1, 57-67. \vspace{0.3cm}

\noindent 6. Grinenko M.M., Birational automorphisms of a
three-dimensional double cone. Sbornik: Mathematics. {\bf 189}
(1998), no. 7, 37-52. \vspace{0.3cm}

\noindent 7. Grinenko M.M., Birational properties of pencils of
del Pezzo surfaces of degrees 1 and 2. Sbornik: Mathematics. {\bf
191} (2000), no. 5, 17-38. \vspace{0.3cm}

\noindent 8. Grinenko M.M., On del Pezzo fibrations. Mathematical
Notes. {\bf 69} (2001), no. 4, 550-565. \vspace{0.3cm}

\noindent 9. Iskovskikh V.A., Rational surfaces with a pencil of
rational curves. Matem. Sbornik. 1967. V. 74 (116), 608-638
(Russian), Engl. transl. in: Math. USSR-Sbornik, {\bf 3} (1967).
\vspace{0.3cm}

\noindent 10. Iskovskikh V.A., Rational surfaces with a pencil of
rational curves and with a positive square of canonical class.
Matem. Sbornik. 1970. V. 83 (125), 90-119 (Russian), Engl. transl.
in: Math. USSR-Sbornik, {\bf 12} (1970). \vspace{0.3cm}

\noindent 11. Iskovskikh V.A., Birational automorphisms of
three-dimensional algebraic varieties, J. Soviet Math. {\bf 13}
(1980), 815-868. \vspace{0.3cm}

\noindent 12. Iskovskikh V.A., On the rationality problem for
three-dimensional algebraic varieties fibered into del Pezzo
surfaces. Proc. Steklov Inst. {\bf 208} (1995), 128-138.
\vspace{0.3cm}

\noindent 13. Iskovskikh V.A., On the rationality criterion for
conic bundles. Sbornik: Mathematics. {\bf 187} (1996), no. 7,
75-92. \vspace{0.3cm}

\noindent 14. Iskovskikh V.A. and Manin Yu.I., Three-dimensional
quartics and counterexamples to the L\" uroth problem, Math. USSR
Sb. {\bf 86} (1971), no. 1, 140-166. \vspace{0.3cm}

\noindent 15. Manin Yu. I. Rational surfaces over perfect fields.
II.  Mat. Sb. {\bf 72} (1967), 161-192. \vspace{0.3cm}

\noindent 16. Manin Yu. I., Cubic forms. Algebra, geometry,
arithmetic. Second edition. North-Holland Mathematical Library,
{\bf 4.} North-Holland Publishing Co., Amsterdam, 1986.
\vspace{0.3cm}

\noindent 17. Pukhlikov A.V., Birational automorphisms of
three-dimensional algebraic varieties with a pencil of del Pezzo
surfaces, Izvestiya: Mathematics {\bf 62}:1 (1998), 115-155.
\vspace{0.3cm}

\noindent 18. Pukhlikov A.V., Birational automorphisms of Fano
hypersurfaces, Invent. Math. {\bf 134} (1998), no. 2, 401-426.
\vspace{0.3cm}

\noindent 19. Pukhlikov A.V., Birationally rigid Fano double
hypersurfaces, Sbornik: Mathematics {\bf 191} (2000), No. 6,
101-126. \vspace{0.3cm}

\noindent 20. Pukhlikov A.V., Birationally rigid Fano fibrations,
Izvestiya: Mathematics {\bf 64} (2000), 131-150. \vspace{0.3cm}

\noindent 21. Pukhlikov A.V., Certain examples of birationally
rigid varieties with a pencil of double quadrics. Journal of Math.
Sciences. 1999. V. 94, no. 1, 986-995.
 \vspace{0.3cm}

\noindent 22. Pukhlikov A.V., Birational automorphisms of
algebraic varieties with a pencil of double quadrics. Mathematical
Notes. {\bf 67} (2000), 241-249. \vspace{0.3cm}

\noindent 23. Pukhlikov A.V., Essentials of the method of maximal
singularities, in ``Explicit Birational Geometry of Threefolds'',
London Mathematical Society Lecture Note Series {\bf 281} (2000),
Cambridge University Press, 73-100. \vspace{0.3cm}

\noindent 24. Pukhlikov A.V., Birationally rigid Fano complete
intersections, Crelle J. f\" ur die reine und angew. Math. {\bf
541} (2001), 55-79. \vspace{0.3cm}

\noindent 25. Pukhlikov A.V., Birationally rigid iterated Fano
double covers. Izvestiya: Mathematics. {\bf 67} (2003), no. 3,
555-596. \vspace{0.3cm}

\noindent 26. Sarkisov V.G., Birational automorphisms of conic
bundles, Izv. Akad. Nauk SSSR, Ser. Mat. {\bf 44} (1980), no. 4,
918-945 (English translation: Math. USSR Izv. {\bf 17} (1981),
177-202). \vspace{0.3cm}

\noindent 27. Sarkisov V.G., On conic bundle structures, Izv.
Akad. Nauk SSSR, Ser. Mat. {\bf 46} (1982), no. 2, 371-408
(English translation: Math. USSR Izv. {\bf 20} (1982), no. 2,
354-390). \vspace{0.3cm}

\noindent 28. Sarkisov V.G., Birational maps of standard ${\mathbb
Q}$-Fano fibrations, Preprint, Kurchatov Institute of Atomic
Energy, 1989. \vspace{0.3cm}

\noindent 29. Sobolev I. V., On a series of birationally rigid
varieties with a pencil of Fano hypersurfaces. Mat. Sb. {\bf 192}
(2001), no. 10, 123-130 (English translation in Sbornik: Math.
{\bf 192} (2001), no. 9-10, 1543-1551). \vspace{0.3cm}

\noindent 30. Sobolev I. V., Birational automorphisms of a class
of varieties fibered into cubic surfaces. Izv. Ross. Akad. Nauk
Ser. Mat. {\bf 66} (2002), no. 1, 203-224.
\vspace{0.3cm} }

\end{document}